\documentclass[11pt,reqno]{amsart}
\usepackage{enumerate}
\usepackage{amsfonts}
\usepackage{amsmath}
\usepackage{amsthm}
\usepackage{amssymb}
\usepackage{latexsym}
\usepackage{multicol}
\usepackage{verbatim}
\usepackage{amscd}
\usepackage{hyperref}
\usepackage{pb-diagram}

\advance\textwidth by 1.2in \advance\oddsidemargin by -.6in \advance\evensidemargin by -.6in
\newtheorem*{cor}{Corollary}
\newtheorem*{lem}{Lemma}
\newtheorem*{prop}{Proposition}

\theoremstyle{definition}
\newtheorem*{defn}{Definition}
\theoremstyle{definition}
\newtheorem{thm}{Theorem}

\newtheorem*{rem}{Remark}

\newenvironment{pf}{\proof}{\endproof}
\newcounter{cnt}
\newenvironment{enumerit}{\begin{list}{{\hfill\rm(\roman{cnt})\hfill}}{%
\settowidth{\labelwidth}{{\rm(iv)}}\leftmargin=\labelwidth%
\advance\leftmargin by \labelsep\rightmargin=0pt\usecounter{cnt}}}{\end{list}} \makeatletter
\def\mydggeometry{\makeatletter\dg@YGRID=1\dg@XGRID=20\unitlength=0.003pt\makeatother}
\makeatother \theoremstyle{remark}

\numberwithin{equation}{section}

 \DeclareMathOperator{\Ht}{ht} 
\let\bwdg\bigwedge
\def\bigwedge{{\textstyle\bwdg}}

\newcommand{\wt}{\operatorname{wt}}

\newcommand{\nc}{\newcommand}
\newcommand{\rnc}{\renewcommand}

\nc{\cal}{\mathcal} \nc{\goth}{\mathfrak} \rnc{\bold}{\mathbf}

\newcommand{\supp}{\operatorname{supp}}

\nc\bomega{{\mbox{\boldmath $\omega$}}} \nc\bpsi{{\mbox{\boldmath $\Psi$}}}
 \nc\balpha{{\mbox{\boldmath $\alpha$}}}
 \nc\bpi{{\mbox{\boldmath $\pi$}}}

\newcommand{\lie}[1]{\mathfrak{#1}}
\makeatletter
\def\section{\def\@secnumfont{\mdseries}\@startsection{section}{1}%
  \z@{.7\linespacing\@plus\linespacing}{.5\linespacing}%
  {\normalfont\scshape\centering}}
\def\subsection{\def\@secnumfont{\bfseries}\@startsection{subsection}{2}%
  {\parindent}{.5\linespacing\@plus.7\linespacing}{-.5em}%
  {\normalfont\bfseries}}
\makeatother

 \nc{\Hom}{\operatorname{Hom}}
\nc{\End}{\operatorname{End}} \nc{\wh}[1]{\widehat{#1}} \nc{\Ext}{\operatorname{Ext}} \nc{\ch}{\text{ch}} \nc{\ev}{\operatorname{ev}}
\nc{\Ob}{\operatorname{Ob}} \nc{\soc}{\operatorname{soc}} \nc{\rad}{\operatorname{rad}} \nc{\head}{\operatorname{head}}

\def\mult

 \nc\boa{\bold a} \nc\bob{\bold b} \nc\boc{\bold c} \nc\bod{\bold d} \nc\boe{\bold e} \nc\bof{\bold f} \nc\bog{\bold g}
\nc\boh{\bold h} \nc\boi{\bold i} \nc\boj{\bold j} \nc\bok{\bold k} \nc\bol{\bold l} \nc\bom{\bold m} \nc\bon{\bold n} \nc\boo{\bold o}
\nc\bop{\bold p} \nc\boq{\bold q} \nc\bor{\bold r} \nc\bos{\bold s} \nc\bou{\bold u} \nc\bov{\bold v} \nc\bow{\bold w} \nc\boz{\bold z}
\nc\boy{\bold y} \nc\ba{\bold A} \nc\bb{\bold B} \nc\bc{\bold C} \nc\bd{\bold D} \nc\be{\bold E} \nc\bg{\bold G} \nc\bh{\bold H} \nc\bi{\bold I}
\nc\bj{\bold J} \nc\bk{\bold K} \nc\bl{\bold L} \nc\bm{\bold M} \nc\bn{\bold N} \nc\bo{\bold O} \nc\bp{\bold P} \nc\bq{\bold Q} \nc\br{\bold R}
\nc\bs{\bold S} \nc\bt{\bold T} \nc\bu{\bold U} \nc\bv{\bold V} \nc\bw{\bold W} \nc\bz{\bold Z} \nc\bx{\bold x}

\newcommand{\mk}{\medskip}

\newcommand{\ZZ}{\bz}

\newcommand{\CC}{\bc}

\newcommand{\yl}{15pt}

\newcommand{\ffbox}[1]{
\setbox9=\hbox{$\scriptstyle\overline{1}$}
\framebox[\yl][c]{\rule{0mm}{\ht9}${\scriptstyle #1}$}
}

\newcommand{\Glie}{\mathfrak{g}}

\newcommand{\Yim}{\mathcal{Y}}

\newcommand{\U}{\cal{U}}
\nc\bvpi{{\mbox{\boldmath $\varpi$}}}

\begin{document}
\setcounter{tocdepth}{1}
\title{Beyond Kirillov--Reshetikhin modules}

\author{Vyjayanthi Chari}
\address{Department of Mathematics, University of California, Riverside, CA 92521, USA}
\email{chari @ math . ucr . edu}
\urladdr{http://math.ucr.edu/\textasciitilde chari}
\author{David Hernandez}
\address{CNRS - Ecole Normale Sup\'erieure, 45 rue d'Ulm, 75005 Paris, FRANCE}
\email{David.Hernandez @ ens . fr}
\urladdr{http://www.dma.ens.fr/\textasciitilde dhernand}

\begin{abstract}

In this survey, we shall be concerned with the category of finite--dimensional representations of the  untwisted quantum affine algebra when the quantum parameter $q$ is not a root of unity. We review  the foundational results of the subject, including the Drinfeld presentation, the classification of simple modules and $q$-characters. We then concentrate on particular families of irreducible representations
 whose structure has recently been understood: Kirillov-Reshetikhin modules, minimal affinizations and beyond.

\vskip 4.5mm

\noindent {\bf 2000 Mathematics Subject Classification:} Primary 17B37, Secondary 81R50, 82B23, 17B67.

\end{abstract}

\maketitle

\tableofcontents
\section{Introduction}

The representation theory of quantum affine algebras has been an active area of research for the past fifteen or twenty years. As in the case of the affine Kac--Moody algebras, there are two distinct but equally important families of representations that are studied: the positive level representations and the level zero representations. The positive level representations are usually studied via the Chevalley generators and Serre relations (in the quantum case this is known as the Drinfeld--Jimbo presentation) while the level zero representations are studied via the loop realization (Drinfeld realization) of the affine (quantum affine) Kac--Moody algebra.

In this survey, we shall be concerned with a full subcategory of the level zero representations of the quantum affine algebra: namely, the category of finite--dimensional representations and we shall assume moreover that the quantum parameter $q$ is not a root of unity. The non--quantum version of this category had been studied previously  in \cite{Cinv, CP87, R} for instance. In those papers, the irreducible finite--dimensional representations of the affine algebra were classified.  An explicit description of the irreducible modules were given in terms of the underlying simple Lie algebra thus
 allowing one to deduce character formulas for the irreducible representations.  The results of those papers, after a straightforward reformulation, prove that the irreducible representations are given by an $n$--tuple of polynomials in an indeterminate $u$, where $n$ is the rank of the underlying simple Lie algebra. It was also clear that the category of finite--dimensional representations was not semi-simple and that there should be some rich structure theory analogous to the representation theory of algebraic groups in characteristic $p$. The blocks of the category were determined in \cite{CM}.

After the presentation of quantum affine algebras was given in  \cite{Dri2}, an identical  classification of the simple finite--dimensional modules for quantum affine algebras was given in \cite{Cha0, Cha}. In the case of the quantum affine algebra corresponding to $\lie{sl}_2$, the simple modules could be described explicitly in a manner similar to the one for the non--quantum case. Outside this case however, it was known that such a result would not hold in general; for instance in \cite{Dri2}, an example was given (in the case of the closely related Yangians) of an irreducible representation of the simple Lie algebra  which did not admit the structure of a module for the quantum affine algebra. In \cite{Cha}, the structure (regarded as a module for the finite--dimensional  simple algebra) of the irreducible \lq\lq fundamental representations\rq\rq of the quantum affine algebra  was determined. In particular, outside the case of $\lie{ sl}_{n+1}$ most of these representations were highly reducible for the simple Lie algebra. It was clear that in the quantum situation  that the structure and characters of the irreducible representations  was much more complex and thus of intrinsic mathematical interest. Together with this, was the external motivation coming from the work of A.N. Kirillov and  N.Reshetikhin, \cite{kr}, the work  of I. Frenkel and N. Reshetikhin \cite{FR92} and others in mathematical physics. The results of these papers suggested that it would be fruitful to concentrate on understanding particular families of irreducible representations.  Thus, the work of \cite{kr} was connected with the irreducible representations of quantum affine algebras corresponding to a multiple of a fundamental representation. These modules are now called the Kirillov--Reshetikhin modules.

The work of \cite{FR92} suggested that there should be some natural minimal family of modules for the quantum affine algebra. This can be also seen in another, purely representation theoretic way: what is the \lq\lq smallest\rq\rq    representation of the quantum affine algebra corresponding to a given  irreducible representation of the simple Lie algebra. This motivated the first author to introduce in \cite{Chari2} the notion of a minimal affinization and these were further studied in \cite{Cha7, Cha8, Cha9}. To do this, one introduces a poset for each dominant integral weight, such that each element of the poset determines a family of  irreducible representation of the quantum affine algebra. A minimal affinization is one which corresponds to the minimal element of the poset.  Another very general family of representations, the so--called prime representations was introduced in \cite{cpp}. The definition is very natural in the light of the results in \cite{Cha0}. It was clear from the results of \cite{Cha7, Cha8} that the Kirillov--Reshetikhin modules are prime and that they  are the  minimal affinizations of multiples of the fundamental weights. 
We shall now see how this and other results motivates the title of our survey.

A very important advance in the theory of finite--dimensional representations was the definition of $q$--characters introduced in \cite{Fre} and studied further in \cite{Fre2, Fre3}. The $q$--characters are generalizations of the usual character of a representation of a simple Lie algebra and satisfy all the usual nice properties: they are additive on direct sums and multiplicative on tensor products. To define them, one observes that in the affine case there is  an infinite family of commuting elements and the $q$--characters encode the information of the generalized eigenvalues and eigenspaces for their action. It was shown in those papers that the irreducible representations were determined by their $q$--character and also that the eigenvalues are given by $n$--tuples of rational functions in an indeterminate $u$.
However, closed formulas such  as an analog of   the Weyl--character formula are not known for $q$--characters. Instead, in \cite{Fre2} an algorithm (now called the Frenkel--Mukhin algorithm)  was proposed to calculate the $q$--character. They proved that the algorithm worked for the fundamental representations and it was thought that this procedure might work in general. This is now known to be false through the recent work of \cite{nn4}. However, it was quite reasonable to expect that there were families of representations for which the algorithm worked and this leads to the notion of a regular representation: one where the Frenkel--Mukhin algorithm yields the $q$--character of the module. It was proved in \cite{Fre2} that a representation was regular if it was special, i.e had an unique eigenvalue given by $n$--tuple of polynomials. From a Lie--theoretic point of view it is more natural to call these minuscule representations as we explain in Section 4.  In \cite{Nab, Nad} it was proved for simply--laced algebras, by using geometric methods, that the Kirillov--Reshetikhin modules were minuscule. An algebraic proof was given in \cite{her06} for all simple Lie algebras.
These results imply the  Kirillov--Reshetikhin conjecture, which gives a closed formula for the character of the tensor product of  Kirillov--Reshetikhin modules. Later in \cite{herma} it was shown that most minimal affinizations are minuscule.

Another way to study the simple finite--dimensional modules for a quantum affine algebras is via branching rules: namely determining the multiplicity of an  irreducible module for the simple Lie algebra in a given module. In the case of  $\lie {sl}_2$ and for fundamental representations, this has been discussed earlier in the introduction. In the case of algebras of type $C_2$, this was done for minimal affinizations  in \cite{Cha2}. For algebras of type $A_n$, it is known that the minimal affinization is irreducible for the simple Lie algebra. In the case of the Kirillov--Reshetikhin modules this problem was studied in \cite{c0} in the case corresponding to nodes of classical type confirming conjectures of \cite{kr, kl, hkoty}. The methods of this paper used the ideas developed in \cite{CPweyl, CPweylq} where the notion of a Weyl module was introduced. The Weyl modules are also  parametrized by an $n$--tuple of polynomials but are not necessarily irreducible but do have very nice universal properties. It was conjectured in \cite{CPweyl, CPweylq} that the Weyl modules are just a tensor product of fundamental representations and we shall see that this is in fact the case in Section 5. 

Relatively little is known about prime representations, in Section 7, we show that most minimal affinizations are prime. One of the main results of \cite{Cbraid} gives a necessary condition for a representation to be prime. Additional examples of prime representations are given in \cite{hle}.
We also discuss  the quasi--minuscule representations, one in which the generalized eigenspaces are of dimension at most one. For some algebras, the Kirillov--Reshetikhin modules and more generally, any minimal affinization is quasi--minuscule. A final family are the small representations closely related
to the geometric small property (Borho-MacPherson). We discuss the
proof \cite{small} of a related conjecture \cite{Nab}
implying a description of the singularities in terms of intersection homology of certain projective morphisms
of quiver varieties.

The study of Kirillov--Reshetikhin modules has been of immense interest in recent years. Character formulas for these representation have been conjectured from physical considerations. These representations are related to several geometric constructions  and to rich combinatorial structures such as crystals \cite{kas, os}, T--systems \cite{kns, hkoty, hkott} and their solutions \cite{Nad, her06} (and more recently cluster algebras \cite{ke, hle}). In the last section we have  a brief overview of these results.

To keep the survey of a manageable length and as a result of our own perspective, we have made many choices and have not elaborated on other important approaches both geometric and combinatorial. In simply-laced cases there are the powerful geometric methods of \cite{Naams}, further studied in \cite{vv}. In particular Nakajima defined an algorithm and proved that it gives the $q$-character of an arbitrary simple representation of a simply-laced quantum affine algebra \cite{Nab}
(a conjectural algorithm for all cases is defined in \cite{her02}). There is also the connection with the theory of crystal basis and details of this approach can be found in the work of \cite{AK, kas, os}. Various important historical references can also be found in \cite{knt}.
We have also restricted ourselves to the case of untwisted affine algebras and when $q$ is not a root of unity. 
This is primarily because the subject is most well--developed in these cases and several results of this survey have not been established for twisted cases. However there are some results in the twisted case which can be found in  \cite{Cha5, her09}. The papers \cite{CProots, bk, Fre4, Nab, h}  consider the case when $q$ is a root of unity. Finally, for the quantum toroidal (ie double affine) case see \cite{her10} and references therein.

\section{Quantum Affine Algebras: definitions and basic results.}\label{un}

Throughout the paper $\bc$ (resp. $\bz$, $\bz_+$, $\bn$) denotes
the set  complex numbers (resp. integers,  non--negative integers,
positive integers) and $q\in\bc$ a fixed non-zero complex number
which is  not a root of unity.

\subsection{} Set $I=\{1,\cdots, n\}$ and $\wh I= I\sqcup\{0\}$. Let $A=(a_{ij})_{i,j\in I}$ be an indecomposable  Cartan matrix of finite type
 and  let  $\hat A=(a_{ij})_{i,j\in \hat I}$   the corresponding untwisted affine Cartan matrix. Fix  a set $\{d_i\}_{i\in\hat I}$ of positive integers so that the matrix $\{d_ia_{ij}\}_{i,j\in \wh I}$ is symmetric.

   Let $\lie g$ and $\wh{\lie g}$ the corresponding finite--dimensional
   simple and untwisted affine Lie algebra associated  to $A$ and $\hat A$ respectively.
As usual,  $R$ denotes the  set of   roots  of $\lie g$ with
respect to a fixed Cartan subalgebra,  and $\{\alpha_i\}_{i\in I}$
(resp. $\{\omega_i\}_{i\in I}$)  a set of simple roots (resp.
fundamental weights). Let $Q$ (resp. $Q^+$) and $P$ (resp. $P^+$)
be the $\bz$--span (resp. $\bz_+$--span) of the simple roots and
fundamental weights respectively and set $R^+= R\cap Q^+$.  Let
$\theta$ be the highest root in $R^+$. Let $\le $ be the usual
partial order on $P$ defined by; for $\lambda,\mu\in P$, we have
$\lambda\le \mu$ iff $\mu-\lambda\in Q^+$.

 Let $W$ be the Weyl group of $R$ and for $i\in I$, let  $s_i\in W$ be the   simple reflection corresponding to
 $\alpha_i$ and let $\bol: W\to \bz_+$ be the length function.
 The  group $W$ acts on the root lattice $Q$ by
  extending $s_i(\alpha_j) = \alpha_j-a_{ij}\alpha_i$.   The affine Weyl group  $\hat{W}$ is isomorphic to the semi-direct product $W\ltimes Q$, under the map $$s_i\to (s_i, 0),\ \ i\in I,\ \
  s_0\to(s_\theta,\theta),$$ where $s_\theta(\alpha_j)=\alpha_j+ a_{0\ j}\theta$. The extended Weyl group $\tilde{W}$ is defined to be the
  semi-direct product $W{{\ltimes}} P$. Regard $W$ and $P$ as subgroups of $\tilde W$ via the maps $w\to (w,0)$ and $\lambda\to t_\lambda=(e,\lambda)$.
  The affine Weyl
  group $\hat{W}$ is a normal subgroup of $\tilde{W}$, and the quotient $\cal{T}=\tilde{W}/\hat{W}$ is a finite group
  isomorphic to a subgroup of the group of diagram automorphisms
  of $\wh{\lie g}$, i.e. the bijections $\tau:\hat I\to\hat I$ such that
  $a_{\tau(i)\tau(j)}=a_{ij}$ for all $i,j\in\hat I$. Moreover, there is
  an isomorphism of groups $\tilde W\cong\cal{T}{\ltimes}\hat W$,
  where the semi-direct product is defined using the action of $\cal{T}$
  in $\hat{W}$ given by $\tau.s_i=s_{\tau(i)}\tau$ (see \cite{Bo}).

\subsection{} For  $i\in\hat I$, set $q_i=q^{d_i}$. For $\ell\in\bz$, and  $r,p,m\in\bz_+$ with $m\ge p$,
set
 $$[\ell]_i=\frac{q_i^\ell-q_i^{-\ell}}{q_i-q_i^{-1}} \text{ , }\ \ [r]_i!=[r]_i[r-1]_i\cdots [1]_i\text{ ,
}\ \ \begin{bmatrix}m\\
p\end{bmatrix}_i=\frac{[m]_i!}{[m-p]_i![p]_i!}.$$ The quantum
affine algebra $\wh{\bu}_q(\lie g)$ is the associative algebra
defined over $\bc$ with generators $k_i^{\pm 1}$,  $x_i^{\pm}$
($i\in \hat I$) and relations:
$$k_ik_j=k_jk_i\text{ , }\ \  k_ix_j^{\pm}=q_i^{\pm a_{i
j}}x_j^{\pm}k_i,$$
$$[x_i^+,x_j^-]=\delta_{i,j}\frac{k_i-k_i^{-1}}{q_i-q_i^{-1}},$$
$${\sum}_{r=0}^{  1-a_{ij}}(-1)^r\begin{bmatrix}1-a_{ij}\\r\end{bmatrix}_{i}
(x_i^{\pm})^{1-a_{ij}-r}x_j^{\pm}(x_i^{\pm})^r=0 \text{ (for
$i\neq j$)}.$$ The
 assignments, \begin{gather*}
\Delta(x_i^+) = x_i^+\otimes k_i + 1\otimes x_i^+,\ \
\Delta(x_i^-) = x_i^-\otimes 1 + k_i^{-1}\otimes x_i^-,\ \
\Delta(k_i^{\pm 1}) = k_i^{\pm 1}\otimes k_i^{\pm 1},\\
S(x_i^+)=-x_i^+k_i^{-1},\ \ S(x_i^-)=-k_ix_i^-,\ \ S(k_i^{\pm 1})= k_i^{\pm 1},\\
\varepsilon (x_i^\pm)=0,\ \ \varepsilon(k_i^{\pm 1})= 1,
\end{gather*}
for $i\in\hat I$, define a Hopf algebra structure on
$\wh{\bu}_q(\lie g)$. The subalgebra  $\bu_q(\lie g)$ generated by
the elements $x_i^\pm$, $k_i$, $i\in I$ is a Hopf subalgebra of
$\wh{\bu}_q(\lie g)$ and is isomorphic as a Hopf algebra to the
quantized enveloping algebra   of $\lie g$. If $\theta=\sum_{i\in
I}\theta_i\alpha_i$, the element $C=k_0\prod_{i\in I}k_i^{-\theta_i}$ is
central in $\wh{\bu}_q(\lie g)$. 

\subsection{}  For $i\in\hat I$ and $m\geq 1$, set $(x_i^\pm)^{(m)} = (x_i^\pm)^m/[m]_i!$
and $(x_i^\pm)^{(0)} = x_i^\pm$.
For $i\in\hat I$ let
 $T_i$ be the algebra automorphism  of $\wh{\bu}_q(\lie g)$
 defined in \cite{L2} by,
  \begin{align*}\label{braid}
    T_i((x_i^+)^{(m)})&=(-1)^mq^{-m(m-1)}(x^-_{i})^{(m)}k_i^m,\ \
    T_i((x^-_{i})^{(m)})=(-1)^mq^{m(m-1)}k_i^{-m}(x^+_{i})^{(m)},\\
    T_i((x_{j}^+)^{(m)})&=\sum_{r=0}^{-ma_{ij}}(-1)^{r}q^{-r}
   ( x_{i}^+)^{(-ma_{ij}-r)}(x_j^+)^{(m)}(x_i^+)^{(r)}\
    \ \text{if $i\ne j$},\\
    T_i((x_j^-)^{(m)})&=\sum_{r=0}^{-ma_{ij}}(-1)^{r}q^{r}
    (x_i^-)^{(r)}
    (x_j^-)^{(m)}(x_i^-)^{(-ma_{ij}-r)}\ \ {\text{if
        $i\ne j$}},
\end{align*} for all $m\ge 0$.
The finite group $\cal T$ acts as Hopf algebra automorphisms of
$\wh{\bu}_q(\lie g)$ by
$$\tau(x_i^\pm)   = x^\pm_{\tau({i})}\ \
\tau (k_i) =k_{\tau(i)}, \ \ \ \text{for all $i\in\hat I$}.$$
Given any $w\in\tilde W$, write $w=\tau s_{i_1}\cdots s_{i_p}$ where $i_r\in\hat I$
for $1\le r\le p$ and  $s_{i_1}\cdots s_{i_p}$
is reduced and set $T_w=\tau T_{i_1}\ldots T_{i_m}$. Then, $T_w$ is an automorphism of $\wh{\bu}_q(\lie g)$
and depends only on $w$.

\subsection{}\label{refd} Following \cite{B1}, define for $i\in I$, $r\in\bz$,
  elements $x^\pm_{i,r}\in\wh{\bu}_q(\lie g)$, by $$x^\pm_{i,r}= o(i)^rT_{\omega_i}^{\mp r}x_i^\pm ,$$
where $o:I\to\{\pm 1\}$ is a map such that $o(i)=-o(j)$
 whenever $a_{ij}<0$ (it is clear that there are exactly two possible choices for $o$).
 Note that $x_{i,0}^\pm=x_i^\pm$. It is proved in \cite{B1} that  the elements $x^\pm_{i,r}$, $k_i$, $C$, $i\in I$, $r\in\bz$ generate $\wh{\bu}_q(\lie g)$. A precise set of defining relations  in terms of these generators (now called the Drinfeld generators and relations) is given in \cite{Dri2} and proved in \cite{B1}.

  For the purposes of this note we identify the following crucial relations. For $i\in I$, $r\in\bz_+$ and $r\ne 0$,  set $$\psi_{i, \pm r}=(q_i-q_i^{-1}) C^{\mp r/2}[x_{i, \pm r}^+, x_{i,0}^-]. $$  For all $i,j\in I$ and $r,s\in\bz$, we have $$[\psi_{i,r},\psi_{j,s}]\in (C-1)\wh{\bu}_q(\lie g).$$ Let $\wh{\bu}_q^\pm(\lie g)$ (resp.
$\wh{\bu}_q(\lie h)$) be  the subalgebra generated by the
$x_{i,r}^\pm$ (resp. the $k_i^{\pm 1}$, $C^{\pm 1}$, $\psi_{i,
r}$). We have,
\begin{equation}\label{dtriang} \wh{\bu}_q(\lie g)=
\wh{\bu}_q^-(\lie g)\wh{\bu}_q(\lie h)\wh{\bu}_q^+(\lie
g),\end{equation}
We shall also use:
\begin{lem}
The assignment $x^\pm_{i,r}\to x^\mp_{i,-r}$, $\psi_{i,
r}\to\psi_{i,-r}$, $q\to q^{-1}$ extends to a $\bc$--linear
algebra antiautomorphism $\Omega$ of $\wh{\bu}_q(\lie
g)$.\hfill\qedsymbol\end{lem}

  For each $r\in\bz$ and $i\in I$ let  $\bu_{i,r}$ be  the subalgebra of $\wh{\bu}_q(\lie g)$ generated by the elements $x^\pm_{i,r}$, $k_i^{\pm  1}$ and let  $\wh{\bu}_i$ be  the subalgebra generated by $\bu_{i,r}$, $r\in \bz$ and $C^{\pm 1}$. Note that $\wh{\bu}_i$
is isomorphic to $\wh{\bu}_{q_i}(\lie {sl_2})$ and that $\bu_{i,r}$ is isomorphic to $\bu_{q_i}(\lie{sl_2})$ for $r\in\bz$.

\subsection{} Following \cite{CProots}, define for $i\in I$, $r\in\bz_+$ elements $P_{i,\pm r}\in\wh{\bu_q}(\lie g)$ by $P_{i,0}=1$ and
\begin{equation}\label{imagroot} P_{i, r}= -\frac{ k_i^{-1}}{1-q_i^{
-2r}}\sum_{s=1}^r\psi_{i,s} P_{i,r-s},\ \ \  \
P_{i,-r}=\Omega(P_{i,r}),\end{equation} and set
 $$\Psi_i^\pm(u)= k_i^{\pm 1}+\sum_{r\ge 1}\psi_{i,\pm r}u^r,\ \ \ \ \bp_i^\pm (u)=\sum_{r\ge 0}P_{i,\pm r} u^
r,$$ where $u$ is an indeterminate. Let $ X_i^\pm$, $i\in I$,  be
the subspace of  $\wh{\bu}_q(\lie g)$ spanned by the elements
$x_{i,r}^\pm$, $r\in\bz$. The following was proved in
\cite{CProots}.
\begin{prop}\label{pbw}  Let $i\in I$.
\begin{enumerit}
\item[\rm{(i)}] There is an equality of power series $\Psi_i^\pm
(u)\bp_i^\pm(u)=k_i^{\pm 1}\bp_i^\pm(q_i^{\mp 2}u).$
\item[\rm{(ii)}] The subalgebra of $\wh{\bu}_q(\lie g)$ generated
by the elements $\{P_{i,\pm r}: i\in I, r\in\bn\}$   is a
polynomial algebra in these variables.
    \item[\rm{(iii)}] For $i\in I$, $r\in\bz_+$, we have $$(x_i^+)^{(r)}(x_{i,1}^-)^{(r)}= (-1)^r q_i^r k_i^{r}P_{i,r} + X_i^-\wh{\bu}_q(\lie g) X_i^+,$$
     $$(x_{i,-1}^+)^{(r)}(x_{i,0}^-)^{(r+1)}= (-1)^r q_i^{r+1} k_i^r\sum_{s=1}^r  x_{i,s+1}^- P_{i,r-s} + \wh{\bu}_q(\lie g) X_i^+.$$ 
     \hfill\qedsymbol

\end{enumerit}
\end{prop}

\subsection{} In general one does not know explicit formulae for the comultiplication in terms of the generators $x^\pm_{i,r}$. However, the next proposition contains partial informations which is sufficient for our purposes. A proof can be found in \cite{B1, BCP, Da}.
\begin{prop}\label{comult} For $i\in I$, $r\in\bz$ we have,

\begin{gather}\label{comult1} \Delta(x_{i,r}^+)\in \sum_{j\in I}\wh{\bu}_q(\lie g) X_j^+\otimes\wh{\bu}_q(\lie g)+
\wh{\bu}_q(\lie g)\otimes\sum_{j\in I}\wh{\bu}_q(\lie g) X_j^+,\\
 \label{comult2}
\Delta(P_{i,r}) -\sum_{s=0}^r P_{i,r-s}\otimes P_s \in \sum_{j\in
I}\wh{\bu}_q(\lie g)X_j^-\otimes\wh{\bu}_q(\lie g)X_j^+.
\end{gather}\hfill\qedsymbol
\end{prop}

\subsection{} We conclude this section by recalling from \cite{CM} the notion of an $\ell$--weight lattice and an $\ell$--root lattice, and also the definition of a braid group action on these lattices. These are   analogous to the definition of the lattices  $P$ and $Q$ and the Weyl group on action on them. The motivation for these ideas will be clear in the next section.

 Let $u$ be an indeterminate and let $\cal P^+$   be the monoid (under coordinate-wise multiplication) of $I$--tuples of polynomials  in $u$ with coefficients in $\bc$ and constant term one.
 Given $i\in I$, $a\in\bc^*$, let $\bpi_{i,a}\in\cal P^+$ be defined by requiring the $i^{th}$ coordinate to be $1-au$
 and all other coordinates to be one. Clearly, $\cal P^+$ is the free abelian monoid generated by the $\bpi_{i,a}$, $i\in I$, $a\in\bc^*$
  and we let $\cal P$ be the corresponding free abelian group. We call $\cal P$ the $\ell$--weight lattice. Let $\wt:\cal P\to P$ be
  the homomorphism of abelian groups given by setting $\wt\bpi_{i,a}=\omega_i$.
 If $\bvpi=(\bvpi_i)_{i\in I}\in\cal P$, then $\bvpi_i\in\bc(u)$ and can  be written as power series  $$\bvpi_i=1+\sum_{r\in\bn}\bvpi_{i,r}u^r,\ \ \bvpi_{i,r}\in\bc, $$ and we shall use this fact freely without further comment.
   Given $\pi\in\bc[u]$ with constant term one,set  $$\pi^+=\pi, \ \ \pi^-=u^{\text{deg} \pi}\pi(u^{-1})/
(u^{\text{deg} \pi}\pi^+(u^{-1}))|_{u=0},\ \, $$ and if
$\varpi=\pi/\pi'\in\bc(u)$, set $\varpi^\pm=\pi^\pm/(\pi')^\pm$.
 For
$\bvpi=(\varpi_1,\cdots,\varpi_n)\in\cal P$ define elements
$\bvpi^\pm\in\cal P$ by
\begin{equation}\label{bvpipm} \bvpi^\pm=(\varpi_1^\pm,\cdots
,\varpi^\pm_n).\end{equation} The element $\bvpi^-$ should not be
confused with the the elements $(\bvpi)^{-1}$ coming from the group
structure.

\subsection{}  Let $B$ be the braid group associated to $W$. Thus,  $B$ is the group generated by elements
$\tilde{T_i}$ ($i\in I$) and defining relations:
\begin{gather*}
\tilde{T_i}\tilde{T_j} =\tilde{T_j}\tilde{T_i},\ \ \text{if}\ \
a_{ij} =0,\ \
\tilde{T_i}\tilde{T_j}\tilde{T_i} =\tilde{T_j}\tilde{T_i}\tilde{T_j},\ \ \text{if}\ \ a_{ij}a_{ji} =1,\\
(\tilde{T_i}\tilde{T_j})^2= (\tilde{T_j}\tilde{T_i})^2,\ \
\text{if}\ \ a_{ij}a_{ji}=2,\ \ (\tilde{T_i}\tilde{T_j})^3=
(\tilde{T_j}\tilde{T_i})^3,\ \ \text{if}\  \ a_{ij}a_{ji} =3, \ \
i,j\in I.\end{gather*}
 The next proposition is a reformulation of
 \cite[Proposition 3.1]{Cbraid} and can be easily checked.
\begin{prop}\label{braid0} There exists a homomorphism of the  group $B$ to   the automorphism group of $\cal P$ given by:
\begin{gather*}
 (\tilde T_i\bvpi)_i
 =\frac{1}{\varpi_i(q_i^2u)},\quad (\tilde T_i\bvpi)_j  =\varpi_j, {\rm{if}} \ a_{ji}=0,\quad
(\tilde T_i\bvpi)_j =\varpi_j(u)\varpi_i(q_iu),\ \  {\rm{if}}\
a_{ji}=-1,\\ (\tilde T_i\bvpi)_j  =\varpi_j(u)\varpi_i(q^3u)
\varpi_i(qu),\ \ {\rm{if}}\ a_{ji}=-2,\\  (\tilde T_i\bvpi)_j =
\varpi_j(u)\varpi_i(q^5u)\varpi_i(q^3u)\varpi_i(qu),\ \ {\rm{if}}\
a_{ji}=-3,
\end{gather*} where $\bvpi\in\cal P$ and $i,j\in I$.
If  $w\in W$  and $s_{i_1}\cdots s_{i_k}$ is a reduced expression
of $w$, the element $\tilde T_w(\bvpi)= \tilde{T_{i_1}}\cdots
\tilde {T_{i_k}}\bvpi$, is independent of the reduced expression
and we have $$\wt(\tilde {T_w}(\bvpi))=w\wt(\bvpi).$$
 \hfill\qedsymbol
\end{prop}
\subsection{}  For $i\in I$, set $$\balpha_{i,a}=
(\tilde T_{i}(\bpi_{i,a}))^{-1}\bpi_{i,a}.  $$ This is exactly analogous to defining the root $\alpha_i$ by $\alpha_i=\omega_i-s_i\omega_i$, and so we have
$\wt(\balpha_{i,a})=\alpha_i$. In addition
\begin{equation}\label{braidroots} \tilde T_j\balpha_{i,a}=\begin{cases}(\balpha_{i,aq_i^2})^{-1},\ \  j=i,\\ \ \balpha_{i,a},\ \ a_{ij}=0,\\
 \balpha_{i,a}\balpha_{j,aq},\ \ a_{i}=-1, \\  \balpha_{i,a}\balpha_{j,aq}\balpha_{j,aq^3}, \ \ a_{ij}=-2,\\
 \balpha_{i,a}\balpha_{j,aq}\balpha_{j,aq^3},\ \ a_{ij}=-3.
 \end{cases}
\end{equation}
Let $\cal{Q}$ (resp. $\cal Q^+$) be the subgroup (resp. monoid)
generated by $\balpha_{i,a}$, $i\in I$, $a\in\bc^\times$. Set
$\cal Q^-=(\cal Q^+)^{-1}$. The following is now immediate.
\begin{lem}
The action of $B$ on $\cal{P}$ preserves
 $\cal{Q}$.
\hfill\qedsymbol
\end{lem}
Let $\preccurlyeq$ be the partial order  defined on $\cal P$ by:
$\bvpi\preccurlyeq\bvpi'$ iff $\bvpi\in\bvpi'\cal Q^-$.
\begin{rem} The elements $\balpha_{i,a}$ are essentially the
elements $A_{i,a}$ first defined in \cite{Fre} from the quantized Cartan matrix and the analog of the partial order $\preccurlyeq$ also appeared
first in \cite{Fre}. The definitions above were  given in \cite{CM} and turn out to be  a natural Lie--theoretic way to introduce these elements and the
partial order. A presentation of the group $\cal P/\cal Q$ can be found in \cite{CM}. \end{rem}

\subsection{} In the process of writing this note, we found that the proof of one of the statements of Lemma 2.7 in \cite{CM}
 was not complete and we take this opportunity to complete that proof.  \begin{prop} \label{correctcm}
 Let $w\in W$ and $i\in I$ be such that $w\alpha_i\in R^+$. Then $\tilde T_w\balpha_{i,a}\in \cal Q^+$ for all $a\in\bc^\times$.
\end{prop}
\begin{pf}  We first prove the Lemma when $\lie g$ is of rank two. If  $\lie g$ is of type $A_2$,
 using  \eqref{braidroots} we get $$ \tilde T_j\balpha_{i,a}=\balpha_{i,a}\balpha_{j,aq},\ \  \tilde T_i\tilde T_j\balpha_{i,a}=\balpha_{j,aq},\ \ i\ne j.$$
 If $\lie g$ is of type $B_2$ let $\alpha_1$ be the long root and $\alpha_2$ the short root, this time we get \begin{gather*}\tilde T_2\balpha_{1,a}=\balpha_{1,a}\balpha_{2,aq}\balpha_{2,aq^3},\ \ \tilde T_1\tilde T_2\balpha_{1,a}=\balpha_{1,aq^2}\balpha_{2,aq}\balpha_{2,aq^3},\ \
\tilde T_2 \tilde T_1\tilde T_2\balpha_{1,a}=\balpha_{1,aq^2},\\
\tilde T_1\balpha_{2,a}=\balpha_{2,a}\balpha_{1,aq},\ \ \tilde
T_2\tilde T_1\balpha_{2,a}=\balpha_{2,aq^4}\balpha_{1,aq},\ \ \
\tilde T_1\tilde T_2\tilde
T_1\balpha_{2,a}=\balpha_{2,aq^4}.\end{gather*} The case of $G_2$
is similar and we leave the calculation to the reader. This proves
the proposition when $\lie g$ is of rank two. For the general
case, we proceed by induction on the length $\bol(w)$ of $w$. If
$\bol(w)=1$ then $w=s_j$ for some $j\ne i$ and \eqref{braidroots}
shows that induction begins. Let $w=s_{i_1}\cdots s_{i_r}$ be a
reduced expression for $w$ and assume that $w\alpha_i\in R^+$. If
$a_{i_r,i}=0$ then $$\tilde T_w\balpha_{i,a}= \tilde
T_{ws_{i_r}}\balpha_{i,a},$$ and since $\bol(ws_{i_r})<\bol(w)$ we
are done by the induction hypothesis. Otherwise we can  write
$w=w_1w_2$ where $w_2$ is in the group generated by $s_{i_r}, s_i$
and $w_2\alpha_i,\  w_1\alpha_i,\ w_1\alpha_{i_r}\in R^+$. By the
rank two case, we know that $\tilde T_{w_2}\balpha_{i,a}$ is in
the monoid generated by $\balpha_{i,c}$, $\balpha_{i_r,d}$ for
$c,d\in\bc^\times$. Since $w_1\alpha_i,\ w_1\alpha_{i_r}\in R^+$,
it follows again by induction that $\tilde T_{w_1}\tilde
T_{w_2}\balpha_{i,a}\in\cal Q_q^+$ and the proof of the inductive
step is complete.

\end{pf}
\begin{cor} If $\bpi\in\cal P^+$ and $w\in W$, then $\tilde T_w\bpi\preccurlyeq\bpi$.
\end{cor}
\begin{pf} It suffices to prove the corollary when $\bpi=\bpi_{i,a}$ for some $i\in I$, $a\in\bc^\times$.
We proceed by induction on $\bol(w)$. If $w=s_i$ then $\tilde
T_i(\bpi_{i,a})=\bpi_{i,a}\balpha_{i,a}^{-1}$
 by definition and  if $w=s_j$, then $\tilde T_j(\bpi_{i,a})=\bpi_{i,a}$. If $\bol(w)>1$ write $w=w's_{i_r}$ for some
 $w'$ with $\bol(w')<\bol(w)$. If $i_r\ne i$, then $\tilde T_w\bpi_{i,a}=\tilde T_{w'}\bpi_{i,a}$ and the result
  follows by induction. If $i_r=i$ then $w'\alpha_i\in R^+$ and we
  have $$\tilde T_w\bpi_{i,a}= \tilde
  T_{w'}\bpi_{i,a} (\tilde T_{w'}\balpha_{i,a})^{-1},$$
  and the result follows from the inductive hypothesis and the proposition.
\end{pf}
\section{The category $\wh{\mathcal{F}}_q(\lie g)$:   simple modules and Weyl modules}

 Let $\wh{\cal F}_q(\lie g)$ be the category   whose objects are finite--dimensional $\wh{\bu}_q(\lie g)$-modules $V$ satisfying,
 $$V=\oplus_{\mu\in P} V_\mu,\ \   V_\mu=\{v\in V: k_i^{\pm 1}v=q_i^{\pm \mu(i)}v,\ \ i\in\hat I\},\ \ \mu=\sum_{i\in I} \mu(i)\omega_i,$$
and where the morphisms between two objects are maps of $\wh{\bu}_q(\lie
g)$--modules.
 Since $\wh{\bu}_q(\lie g)$ is a Hopf algebra the category
 $\wh{\cal F}_q(\lie g)$ is closed under taking tensor products and duals and we let
 $\wh{\text{Rep}}(\lie g)$ be the corresponding  Grothendieck ring.
   This category
   is far from semi-simple and a parametrization of the blocks of the category can be found in
   \cite{CM}.

    In this section we shall see that given $\bpi\in\cal P^+$
   we can associate to it canonically two  modules: $V(\bpi)$ which is simple
    and $W(\bpi)$ which has nice universal properties and see that
    $\cal P^+$ parametrizes the simple objects in $\wh{\cal
    F}_q(\lie g)$. These results can be found in \cite{Cha0, Cha, CPweyl, CPweylq}.  We include a sketch of a proof in
some cases for the readers convenience and for motivating some of
the later results.

   We shall also begin the discussion  the $q$--character of objects in $\wh{\cal F}_q(\lie g)$.
     These were   originally defined and studied in a slightly different formulation in \cite{Fre}.
 For our purposes, we define it as follows. It is not hard to see \cite{Cha0}, that the  element $C-1$ of $\wh{\bu}_q(\lie g)$ acts trivially on an object of $\wh{\bu}_q(\lie g)$ and hence the elements $k_i^{\pm 1}$, $P_{i,r}$, $i\in I$,
 $r\in\bz$ act as  a family of commuting operators on $V$. Hence we can write $V$ as a direct
  of generalized eigenspaces for their action, i.e we have
   $$V=\oplus_{\bod} V_\bod,\ \ \bod=(d_{i,r})_{i\in I, r\in\bz},\ \  d_{i,r}\in\bc,\ \
    ,$$ where $$V_\bod=\{v\in V: k_i^{\pm 1} v= q^{\pm d_{i,0}}v,\ \ (P_{i,r}-d_{i,r})^{N_{i,r}}v=0, \ {\rm{for\ some}}\
    N_{i,r}\in\bn\}.$$
    The formal sum \begin{equation}\label{qch} \ch_q(V)=\sum_{\bod}\dim V_\bod \bod,\end{equation}
    is called the $q$--character of $V$.

\subsection{} We shall assume that the reader is aware
 that the category $\cal F_q(\lie g)$ of (type $1$) finite--dimensional representations of
  $\bu_q(\lie g)$ is essentially the same as the
  corresponding category for the simple Lie algebra as
long as $q$ is not a root of unity. Thus, the category is
semi-simple, the irreducible representations are parametrized by
dominant integral weights and given $\lambda=\sum_{i\in
I}\lambda(i)\omega_i\in P^+$, we let $V(\lambda)$ be the
irreducible finite--dimensional module generated by an element
$v_\lambda$, with defining relations:
$$k_i^{\pm 1} v_\lambda=q^{\pm\lambda(i)}v_\lambda,\ \ x_i^+v_\lambda=0,\ \ (x_i^-)^{\lambda(i)+1}v_\lambda=0,\ \ i\in I.$$
 In
particular any finite--dimensional representation $V$ of
$\bu_q(\lie g)$ can be written as
$$V=\oplus_{\mu\in P} V_\mu,\ \ V_\mu=\{v\in V: k_i^{\pm 1}v=q^{\pm \mu(i)}v\} $$
and  we let $$\ch(V)=\sum_{\mu\in P}\dim V_\mu e(\mu),$$ be the
element of the group ring $\bz[P]$. Set $$\wt(V)=\{\mu\in P:
V_\mu\ne 0\}.$$
 Any object of $\cal F_q(\lie
g)$ is completely determined, up to isomorphism,  by its
character. Details of all these facts  can be found in any of the
standard books on quantum groups (for example \cite{Cha2}).
 Since any object $V$ of $\wh{\cal F}_q(\lie g)$ can also be regarded
  as an object in $\cal F_q(\lie g)$ the set $\wt (V)$ is defined and  we see immediately, that $\ch_q(V)$  is much finer than the character $\ch(V)$.

\subsection{} We introduce the notion of a highest weight module
adapted to the triangular decomposition given in \eqref{dtriang}.

\begin{defn} We say that a $\wh{\bu}_q(\lie g)$-module $V$ is $\ell$-highest weight with highest weight vector $v$ and highest weight $\bod=\{d_{i,r}\in\bc: i\in I,\ r\in\bz\}$  if $ V=\wh{\bu}_q(\lie g)v $ and $$x_{i,r}^+v =0,\ \ k_i^{\pm 1}v= q_i^{\pm d_{i,0}} v,\ \ P_{i,r} v= d_{i,r} v,\ \ \ \  r\in\bz, r\ne 0.$$
\end{defn}
We begin with the following result which is an
immediate consequence of Proposition \ref{comult}.
\begin{lem} \label{tp}  Let $V_p$ be $\ell$--highest weights $\bod^p$ and highest weight vectors $v_p$, $p=1,2$. The $\wh{\bu}_q(\lie g)$ submodule of $V_1\otimes V_2$ generated by $v_1\otimes v_2$ is also $\ell$--highest weight, with highest weight vector $v_1\otimes v_2$ and $\ell$--highest weight $\bod$ given by $$d_{i,\pm r}=\sum_{s=0}^r d^1_{i, \pm s}d^2_{i, \pm(r-s)},\  \ i\in I, r\in\bz_+. $$
More generally, if $V_p\in\wh{\cal F}_q(\lie g)$ then the
$q$-character of $V_1\otimes V_2$ is the product of the
$q$--characters of $V_1$ and $V_2$. \hfill\qedsymbol
\end{lem}
\subsection{}  The next result gives a necessary condition for an $\ell$-highest weight module to be finite--dimensional.
\begin{lem} Let $V$ be an $\ell$--highest weight module with highest weight  $\{d_{i,r}: i\in I,\ r\in\bz\}$ and highest weight vector $v$.
Then $\dim(V) <\infty$  only if for all $i\in I$,
 we have: \begin{gather} d_{i, 0} = S_i\in\bz_+,\quad \ d_{i,r}= 0,\  \ |r| >S_i, \ \quad
 d_{i,\pm S_i}\ne 0,\ \ \\ \label{pm} 1+ \sum_{r\ge 1} d_{i,-r}u^r = u^{S_i}+ d_{i,S_i}^{-1}\sum_{r\ge 1} d_{i,r} u^{S_i-r}.  \end{gather}
In other words, if we set $\pi_i=\sum_{r\ge 0}d_{i,r}u^r$ and
$\bpi=(\pi_1,\cdots,\pi_n)$, then
 $\bpi\in\cal P^+$ and the highest weight of $V$ is given by $\bpi$, in the sense that $$\sum_{r\ge 0}d_{i,\pm r}u^r=\bpi^\pm(u),\ \ S_i=\deg\pi_i.$$
\end{lem}
\begin{pf} For $i\in I$, $r\in\bz$, regard $V$ as a module for the subalgebra $\bu_{i,r}$
 which  we recall, is isomorphic to $\bu_{q^i}(\lie{sl_2})$.
 Since $V$ is finite--dimensional  there exists $S_i\in\bz_+$ minimal such that  $(x_{i, r}^-)^{S_i+1} v=0$.
  It is  now immediate from the representation theory of $\bu_q(\lie{sl_2})$ that $S_i+1=d_{i,0}+1$ and hence $d_{i,0}=S_i$.
  To prove that $d_{i,r}=0$ if $|r|>S_i$ we  use Proposition \ref{pbw}(iii)
   by noting that the second term on the right hand side is zero on $v$
    and  the left hand side is zero on $v$ if $r>S_i$ by the preceding discussion.
     Finally \eqref{pm} is obtained by using Proposition \ref{pbw}(i).
\end{pf}

\subsection{} We now turn our attention to the converse problem.  Given $\bpi\in\cal P^+$,  the Weyl module $W(\bpi)$ is defined in \cite{CPweyl} as
the $\wh{\bu}_q(\lie g)$--module generated by an element $w_\bpi$
with relations:
$$x_{i,r}^+ w_\bpi=0,\ k_i^{\pm1} w_\bpi=  q^{\pm\deg\pi_i} w_\bpi,\  \  \ \ \bp^\pm_i(u) w_\bpi=\pi^\pm_i(u) w_\bpi,\ \ (x_{i,0}^-)^{\deg(\pi_i)+1} w_\bpi=0,$$ for all $i\in I$ and $r\in\bz$. The following theorem was proved in \cite{CPweyl}, although part (iii) of the theorem was proved much earlier in \cite{Cha0, Cha}. Note that part (ii) shows that Weyl modules are universal
$\ell$--highest weight modules in $\wh{\cal F}_q(\lie g)$.
\begin{thm} \label{weyldefn}
\begin{enumerit}
\item For $\bpi\in\cal P^+$ we have $W(\bpi)\in\Ob\cal F_q(L(\lie
g))$. Moreover, $$\dim(W(\bpi))_{\wt\bpi}=1,\ \
\wt(W(\bpi))\subset\wt\bpi -Q^+. $$ In particular, $W(\bpi)$ is
indecomposable and  has a unique irreducible quotient $V(\bpi)$.
\item Any $\ell$--highest weight module in $\wh{\cal F}_q(\lie g)$
is a quotient of $W(\bpi)$ for some $\bpi\in\cal P^+$. \item  If
$V\in\Ob \wh{\cal F}_q(\lie g))$ and $V$ is simple, then  $V\cong
V(\bpi)$ for some $\bpi\in\cal P^+$.\end{enumerit}\hfill\qedsymbol
\end{thm}
\begin{rem} Lemma \ref{tp} and  the  preceding theorem makes clear  the motivation for defining the monoid $\cal P^+$, it is the analog of the fact that the monoid $P^+$ parametrizes the
irreducible finite--dimensional representations of $\lie g$.  \end{rem}

\begin{defn} For  $i\in I$, $a\in\bc^\times$, the module $V(\bpi_{i,a})$ is called the $i^{th}$--fundamental module with parameter $a$. \end{defn}
The following is a consequence of the Theorem and Lemma \ref{tp} and justifies the definition of the fundamental module. It is analogous to the result that any object in $\cal F_q(\lie g)$ occurs in the tensor product of the fundamental representations $V(\omega_i)$, $i\in I$.
 \begin{cor}\label{weylsimple} Let $\bpi=\prod_{i\in I}\prod_{s=1}^{\deg\pi_i}\bpi_{i,b_{i,s}}\in\cal P^+$. The module $V(\bpi)$  is a subquotient of the tensor product (in arbitrary order) of the modules  $V(\bpi_{i,b_{i,s}})$ $i\in I$, $b_{i,s}\in \bc^\times$, $1\le s\le\deg\pi_i$.\end{cor}

\section{The category $\wh{\cal F}_q(\lie{sl_2})$ and  $q$--characters as elements of $\bz[\cal P]$ } In this section, we shall see that the formal sum
  \begin{equation}\label{qcha} \ch_q(V)=\sum_{\bod}\dim V_\bod \bod,\end{equation} can be regarded
     as an element of the integral group ring of $\cal P$ and
     explain how this is related to the original formulation of
     $q$--characters in
     \cite{Fre}. This requires an understanding of the simple
     objects of $\wh{\cal F}_q(\lie {sl_2})$ and we also discuss the
     $q$--characters of the irreducible modules and the Weyl modules in the $\lie{sl_2}$--case.

  \subsection{} It is convenient
  at this point to introduce some definitions.  Given $V\in\Ob\wh{\cal F_q}(\lie g)$
  and $\bvpi\in\cal P$ recall the elements  $\bvpi^\pm=(\varpi_1^\pm,\cdots \varpi^\pm_n)$ defined in \eqref{bvpipm} and set $$V_\bvpi= V_{\{\varpi^\pm_{i,r}: i\in I, r\in\bz\}},\ \ \varpi_i^\pm(u) =\sum_{r\ge 0}\varpi^\pm_{i,r}u^r,$$
 Note that since $\dim W(\bpi)_{\wt\bpi}=1 =\dim V(\bpi)_{\wt\bpi}$, it follows also that  $\dim W(\bpi)_\bpi =1=\dim V(\bpi)_\bpi$.

 \begin{defn} We say that $V\in\Ob\wh{\cal F_q}(\lie g)$ is {\em
 minuscule}
 if there exists exactly one element
 $\bpi\in\cal P^+$ such that $\dim V_\bpi\ne 0$. We say that $V\in\Ob\wh{\cal F_q}(\lie g)$
 is {\em quasi--minuscule} if for all $\bpi\in\cal P$ we have $\dim V_\bpi\le
 1$. Finally we say that $V\in\Ob\wh{\cal F_q}(\lie g)$ is {\em
 prime}
 iff $V$ cannot be written as a tensor product of nontrivial
 objects of $\wh{\cal F_q}(\lie g)$.
  \end{defn}
 By Theorem \ref{weyldefn}, it is clear that given an $V\in\Ob\wh{\cal F_q}(\lie g)$ there must
 exist $\bpi\in\cal P^+$ such that $V_\bpi\ne 0$. If $V$ is minuscule it follows that any irreducible
  constituent of $V$ must be isomorphic to $V(\bpi)$ and in particular, if $\dim V_\bpi=1$, then $V\cong V(\bpi)$.
   In the literature, so far  simple minuscule modules are called special \cite{Nab} and quasi--minuscule modules are called thin \cite{small}.
 We believe that our notation is more consistent with the representation theory of simple Lie algebras, where a minuscule
    representation is one with a unique dominant integral weight and a quasi--minuscule is one where all weight spaces are of dimension at most one.
    In fact it is quite easy to see that if $V\in\Ob\wh{\cal F}_q(\lie g)$ is minuscule
    (resp. quasi--minuscule) as an object of $\cal F_q(\lie g)$ then it is minuscule (resp. quasi-minuscule) as an object of $\wh{\cal F}_q(\lie g)$.
    The converse is far from true as will be clear from the
 examples of minuscule and quasi-minuscule representations given in the rest of the
 paper. Finally, note that if $V\in\Ob\wh{\cal F}_q(\lie g)$ is
 irreducible when regarded as an object of $\cal F_q(\lie g)$
 then $V$ is prime. The modules $V(\bpi_{i,a})$ for $i\in I$ and
 $a\in\bc^\times$ are clearly prime but as we shall see in general
 not irreducible as an object of $\cal F_q(\lie g)$.

\subsection{}  As in the case of the representation theory of simple Lie algebras
 it is crucial to first understand the case when $\lie g$ is of type $\lie sl_2$.
 In this case an  element of $\cal P^+$ is just a single polynomial $\pi$ in $u$ with
 constant term one. Our goal is to describe the structure of $V(\pi)$ and $W(\pi)$.

 Given $a\in\bc^\times$ and $m\in\bz_+$, set $$\pi^m_ a(u)=\prod_{j=1}^m(1-q^{m-2j+1}au), \ \ m\ge 1, \ \ \pi^0_a=1.$$
 It is proved in \cite{Cha0} that  $V(\pi^m_a)$ is irreducible when regarded as module for
 $\bu_q(\lie{sl_2})$ and we have $$\dim V(\pi^m_a)=m+1,\ \ \wt(V(\pi^m_a))=\{m-2j: 0\le j\le m\}.$$
 Moreover,  if $v\in V(\pi^m_a)_{m-2j}$, we have an equality of power series \cite{CM},
  \begin{gather}\label{qchsl2} P_1^\pm v= ({\pi^{m-j}_{aq^{-j}}})^\pm(({\pi^j_{aq^{m-j+2}}})^\pm)^{-1}v
   = \pi^m_a(\balpha_{aq^{m-1}}\cdots\balpha_{aq^{m-2j+1}})^{-1}v, \end{gather} where $P_1^\pm$ are the elements defined in
  \eqref{imagroot}. The first equality in the preceding formula appears in a different
  form in \cite{Cha0}, while the second appears  in \cite{Fre} using the elements $\psi_{i,\pm r}$ and $A_{i,a}^{-1}$.

   In particular we see that $V(\pi^m_a)$ is minuscule, quasi-minuscule and prime. Moreover we have, $$ \bod=\{d_r\}_{r\in\bz} ,\  V(\pi^m_a)_\bod\ne 0,\ \ \implies \sum_{r\ge 1}d_{\pm r}u^r=\varpi^\pm(u),\ \ {\rm{for\ some}}\ \ \varpi\in\cal P.$$
Suppose now that  $\pi\in\cal P^+$. It is easy to see that  there
exists:
  \begin{enumerit} \item a unique partition of $\deg\pi=(m_1,m_2,\cdots  ,m_s)$ ($m_1\ge m_2\ge\cdots\ge  m_s\ge 1$),
    \item and unique elements $a_1,\cdots, a_s\in\bc^\times$ satisfying  $$a_k/a_p\notin\{ q^{\pm (m_k+m_p)}, q^{\pm (m_k+m_p-2)},\cdots ,q^{\pm (m_k-m_p+2)}\},\ \ p>k,$$
        \end{enumerit} such that $\pi=\prod_{k=1}^s\pi_{m_k}^{a_k}$.
We call this the $q$--factorization of $\pi$ and  the elements
$\pi_{m_k}^{a_k}$ the $q$--factors of $\pi$. The following theorem
was proved in \cite{Cha0}.

 \begin{thm}\label{sl2} Let $\bpi\in\cal P^+$ and assume that $\bpi=\prod_{k=1}^s\pi^{m_k}_{ a_k}$ is a $q$--factorization of $\pi$. Then
$$V(\bpi)\cong V(\pi^{m_1}_{a_1})\otimes \cdots\otimes
V(\pi^{m_s}_{a_s}).$$\end{thm}

\begin{cor} Let $\bpi\in\cal P^+$ and assume that $\bpi=\prod_{k=1}^s\pi^{m_k}_{ a_k}$ is a $q$--factorization of $\pi$. Then,
 $V(\bpi)$ is minuscule iff  for $1\le j\neq k\le m$, we have that,
  $$a_k \notin \{a_jq^{m_j-m_k - 2},\cdots,
a_j q^{- m_j-m_k + 2}\},$$
 and it is quasi-minuscule iff $V(\bpi)$ is minuscule and for $j\neq k$, $$a_k  \neq
a_jq^{m_j-m_k}.$$ Finally the irreducible prime objects in $\wh{\cal
F}_q(\lie g)$ are the $V(\pi_a^m)$, $m\in\bz_+$, $a\in\bc^\times$.
\end{cor}

\begin{pf} The first statement is proved in \cite[Lemma 4]{Fre} while the last  is a direct consequence of Theorem \ref{sl2}. For the second,  suppose that $a_kq^{m_k} \in \{a_jq^{m_j},\cdots, a_j q^{- m_j + 2}\}$, say $a_kq^{m_k - 1} = a_j q^{m_j - 1 - 2r}$ where $r\leq m_j - 1$.
So $\bpi \balpha_{a_jq^{m_j}}^{-1}\balpha_{a_jq^{m_j -
2}}^{-1}\cdots \balpha_{a_jq^{m_j - 2r}}^{-1}$ has multiplicity at
least $2$ in $\ch_q(V(\bpi))$ and so $V(\bpi)$ is not
quasi-minuscule.

Otherwise for $j\neq k$, we have $\{a_k q^{m_k - 1},\cdots
a_kq^{1-m_k}\}\cap \{a_jq^{m_j - 1},\cdots, a_jq^{1 - m_j}\} =
\emptyset$ and the quasi-minuscule property follows.
\end{pf}

 \subsection{}\label{rappel} For general $\lie g$, we have the following
 consequence of Theorem \ref{sl2} which first appeared in
 \cite{Fre}.
 \begin{prop} Let $\lie g$ be simple and  $V\in\wh{\cal F}_q(\lie g)$. Then $V=\oplus_{\varpi\in\cal P} V_\varpi$.
 In particular, ${\rm ch_q}
(V)\in\bz[\cal P]$ and
 we have a ring homomorphism ${\rm ch_q} :\wh{\text{Rep}}(\lie g)\rightarrow \bz[\cal P]$.
 \end{prop}
 \begin{pf} In the case when $\lie g$ is $\lie sl_2$ and $V$ is irreducible the Proposition is immediate from Theorem \ref{sl2}.
 If $V$ is reducible then the result follows again since the eigenspaces for the action of
  $P_{1,r}$ on $V$ are just the sum of the eigenspaces coming from a Jordan--Holder series for $V$.
   For arbitrary $\lie g$, the result follows by regarding $V$ as a module
   for $\wh{\bu}_{i}$, $i\in I$.
  The fact that $\ch_q$ is  a ring homomorphism is immediate from Lemma \ref{tp}.
 \end{pf}

\begin{rem} We now explain the connection with original
formulation of \cite{Fre}.  Let $\Yim$ be the polynomial ring over
the integers in the infinitely many
     variables $Y_{i,a}^{\pm 1}$, $i\in I$, $a\in\bc^\times$
      The assignment  $$\bpi = \prod_{i\in I, a\in\CC^*}\bpi_{i,a}^{r_{i,a}}
      \mapsto m = \prod_{i\in I, a\in\CC^*} Y_{i,a}^{r_{i,a}},$$ is
      clearly an isomorphism of rings and hence
      $\ch_q(V)$ can also be regarded as an element of $\cal Y$,
         For some applications the monomial
         notation is more convenient, and in the following we shall  use  both notations without comment.
          Any object $V$ in $\wh{\cal F}_q(\lie g)$ is obviously a finite--dimensional module for $\bu_q(\lie g)$
and it follows from the fact that $\ch_q(V)\in\bz[\cal P]$ that
$\ch(V)=\wt(\ch_q(V)),$ or equivalently \cite{Fre}, in the monomial notation
one replaces  $Y_{i,a}$ by $y_i = e(\omega_i)$ to get the usual
character.
\end{rem}

 \subsection{} We go back to the case $\lie g$ of type $\lie sl_2$ and turn to understanding the Weyl modules. It was proved in \cite{CPweylq} that if $\pi=\prod_{r=1}^s\pi^1_{a_r}$,  and $a_r/a_{k}\ne q^{-2}$ if $k>r$, then \begin{equation}\label{weylsl2}W(\pi)\cong V(\pi^1_{a_1})\otimes\cdots\otimes V(\pi^1_{a_r}).\end{equation}
We now compute the $q$--character of $W(\pi)$.  Define a map $$
\cal P^+(\pi) =\{\pi'\in\cal P^+: \pi(\pi')^{-1}\in\cal P^+\} \to
\cal P,\ \ \pi'\mapsto \pi'(u)\pi'(q^2u)(\pi(q^2u))^{-1}. $$  We
claim that the map is injective. For, if
$$\pi'(u)\pi'(q^2u)=\eta'(u)\eta'(q^2u),$$ then choose
$a\in\bc^\times$ so that $(1-au)$ divides the left hand side and
$(1-aq^ku)$ does not divide it for any $k<0$. It is immediate that
$(1-au)$ divides $\pi'(u)$ and similarly must divide $\eta'(u)$.
An obvious induction on $\deg \pi'$ now proves that $\pi'=\eta'$.
Given,  $\pi\in\cal P^+$ and $\pi'\in\cal P^+(\pi)$, write
$$\pi=\prod_{j=1}^k(1-a_ju)^{r_j},  \ \
\pi'=\prod_{j=1}^k(1-a_ju)^{S_j},$$ where $a_j\in\bc^\times$,
$1\le j\le k$  are distinct, and set
$$d_{\pi}(\pi')=\binom{r_1}{S_1}\binom{r_2}{S_2}\cdots\binom{r_k}{S_k}.$$
Using \eqref{qchsl2}, \eqref{weylsl2} and the multiplicative
property of $\ch_q$, we see that 
$$\{\varpi\in \cal P:
W(\pi)_\varpi\ne 0\}=\{ \pi'(u)\pi'(q^2u)(\pi(q^2u))^{-1}:
\pi'\in\cal P^+(\pi)\},\ $$   and we have proved the following:
\begin{prop}\label{weylsl2qm} Let $\pi\in\cal P^+$. We have $${\rm{\ch}}_q(W(\pi))=\sum_{\pi'|\pi}d_{\pi}(\pi')\pi'(u)\pi'(q^2u)(\pi(q^2u))^{-1}. $$In particular, $W(\pi)$ is quasi--minuscule iff $\pi$  has distinct roots.
\end{prop}

\section{Irreducibility of Tensor products and $q$--characters of Weyl modules}

In this section, we give a partial analog of Theorem \ref{sl2} for
general $\lie g$. Thus we  give a sufficient condition
\cite{Cbraid}  for a tensor product of two irreducible objects of
$\wh{\cal F}_q(\lie g)$ to be irreducible and hence also a
necessary condition for a representation to be prime. We shall
then see that this result along with results in \cite{kas, Naams, bn} can be used to substantially strengthen
Corollary \ref{weylsimple}: namely that any $\ell$--highest weight
module is
 actually a quotient of an ordered tensor product of fundamental modules.
 This is essentially
 equivalent
  (by the universal property of Weyl modules) to proving
  that for all $\bpi\in \cal P^+$ the module $W(\bpi)$ is   isomorphic to
  an ordered tensor product of fundamental modules.   Finally, we discuss $q$--characters of the
     fundamental representations and hence also the
     $q$--characters of Weyl modules.

\subsection{}  We shall say that two polynomials $\pi$ and $\pi'$ with constant term one are in general position if the $q$--factorization of $\pi\pi'$ is just the product of the $q$--factorization of $\pi$ and $\pi'$.  We say that $\pi$ is in general position with respect to $\pi'$ if for every  $q$--factor $\pi^m_a$ of $\pi$ and $\pi^r_b$ of $\pi'$ we have $$a\ne bq^{-(m+r-2p)},\ \ 0\le p<\min(m,r).$$

The following was proved in \cite{Cbraid}. It was motivated by the following classical result: if $V$ is any finite--dimensional $\lie g$--module and $w\in W$, then $\dim V_\mu=\dim V_{w\mu}$ for all $\mu\in P$. In the case of the quantum affine algebra it was natural to ask if something similar was true for $\wt_\ell(V)$. This motivated the definition of the braid group action on $\cal P$ (in the quantum case, it is quite natural to replace the Weyl group by the braid group). Part (i) of the theorem shows that the analog of the classical result works for the highest weight $\bpi$; the discussion in the $\lie sl_2$ case shows that this is false for an arbitrary weight. However, even this partial information is enough to find conditions for a tensor product of irreducible representations to be irreducible.

\begin{thm}\label{braidthm}  Let $\bpi,\bpi'\in\cal P^+$. For all $w\in W$ we have \begin{equation}\label{braidinv}  \dim V(\bpi)_{T_w\bpi}=\dim V(\bpi)_{\wt(T_w\bpi)}=1.\end{equation}
 Further, $V(\bpi)\otimes V(\bpi')$ is an $\ell$--highest weight module (resp. irreducible) if given a reduced expression  $w_0=s_{i_1}\cdots s_{i_N}$ of the longest element in $W$, we have
\begin{enumerit}\item the elements $\pi_{i_N}$ is in general position with respect to $\pi'_{i_N}$  (resp. $\pi_{i_N}$ and  $\pi'_{i_N}$ are in general position),
\item For $2\le j\le N$, the elements
$(\tilde{T}_{i_j}\cdots\tilde{T}_{i_N}\bpi)_{i_{j-1}}$ is in
general position with respect to $\pi'_{i_{j-1}}$ (resp.
$(\tilde{T}_{i_j}\cdots\tilde{T}_{i_N}\bpi)_{i_{j-1}}$ and
$\pi'_{i_{j-1}}$ are in general position).
\end{enumerit}

\end{thm}
\begin{pf} The idea of the proof is as follows. In the case when $\lie g$ is of type $\lie {sl_2}$, this was done in \cite{Cha0}. Then one observes that for all $w\in W$, the weight space $\dim V(\bpi)_{w\wt\bpi}=1$ and hence is an eigenspace for the action of the  elements $P_{i,r}$, $i\in I$, $r\in\bz_+$. The corresponding   eigenvalue  is $\tilde{T}_w(\bpi)$. Given a reduced expression $s_{i_1}\cdots s_{i_N}$ of $w_0$, set $v_{N}=v_\bpi$ and fix non--zero elements $v_j\in V(\bpi)_{s_{i_{j+1}}\cdots s_{i_N}\wt\bpi}$ for $1\le j\le N-1$. It is easy to see that
 $x_{j,r}^+ v_j=0$ for all $r\in\bz_+$. We proceed by induction to show that $v_j\otimes v_{\bpi'}\in \wh{\bu}_q(\lie g)(v_\bpi\otimes v_{\bpi'})$. The  inductive step follows from the fact that the element $v_j\otimes v_{\bpi'}$ is an $\ell$--highest weight vector for the subalgebra $\wh{\bu}_i$ (which is isomorphic to $\wh{\bu}_{q^i}(\lie {sl}_2)$) and by using condition (ii) in the statement of the theorem.
\end{pf}

Let $\bpi\in\cal P^+$ and write $\bpi=\prod_{i\in
I}\prod_{p=1}^{S_i}\bpi_{i,a_{i,p}}$ for some $a_{i,p}\in\bc^\times$.
Fix an ordering $\le $ of the set $\{a_{i,p}: i\in I, 1\le p\le
S_i\}$ so that the polynomial $(1-a_{i,p}u)$ is in general
position with respect to $(1-a_{j,r}u)$ if $a_{i,p}<a_{j,r}$.
\begin{cor}\label{fundirr}  The ordered tensor product $\otimes_{i\in I}\otimes_{p=1}^{S_i} V(\bpi_{i,a_{i,p}})$ is $\ell$--highest weight and hence a quotient of $W(\bpi)$. In particular if $a_{i,p}=1$ for all $i\in I$ and $1\le p\le s$ the tensor product is irreducible. \end{cor}
The corollary had been proved in \cite{AK} when $\lie g$ is of
type $A_n$ or $C_n$,
 a geometric proof was given for the simply--laced algebras in \cite{vv} and a complete proof was
 given in \cite{kas} by using  crystal bases. In \cite{Naams}, Nakajima introduced the notion of standard
 modules for $\wh{\bu}_q(\lie g)$ through the geometry of quiver varieties and the results of \cite{vv}
 show in the simply-laced case,  that an algebraic definition of the standard module is just the tensor product given in  the corollary.

\subsection{} We can now state the following result which can be viewed as giving a presentation of the standard modules.
\begin{thm}\label{presentweyl} Let $\bpi\in\cal P^+$
 and write $\bpi=\prod_{i\in I}\prod_{p=1}^{S_i}\bpi_{i,a_{i,p}}$ for some $a_{i,p}\in\bz_+$. Fix an ordering $\le $ of the set $\{a_{i,p}: i\in I, 1\le p\le S_i\}$ so that the polynomial $(1-a_{i,p}u)$ is in general position with respect to $(1-a_{j,r})$ if $a_{i,p}<a_{j,r}$.
The ordered tensor product $\otimes_{i\in I}\otimes_{p=1}^{S_i}
V(\bpi_{i,a_{i,p}})$ is isomorphic to $W(\bpi)$. \hfill\qedsymbol
\end{thm}
This theorem was conjectured in \cite{CPweyl}. A proof was given
in the case of $\lie{sl_2}$ in \cite{CPweylq}. For arbitrary
simple Lie algebras this follows  easily from the results of \cite{Naams, bn} and a proof can be found in
\cite{CM}. Another approach which is more algebraic and similar to
the proofs given for the $\lie{sl_2}$ case can be found for $\lie
sl_{n+1}$ in  \cite{CL} and  a proof using Demazure modules can be
found in \cite{FoL} in the general simply--laced  case.

  \subsection{} Given $V\in\Ob\wh{\cal F}_q(\lie g)$, set
  $$\wt_\ell(V)=\{\bvpi\in\cal P^+: V_\bvpi\ne 0\}.$$
\begin{thm}\label{lcone} \begin{enumerit} \item For  $i\in I$ and
$a\in\bc^\times$, we have $\wt_\ell
V(\bpi_{i,a})\subset\bpi_{i,a}\cal Q^-.$
 \item For   $\bpi\in\cal P^+$, we have $\wt_\ell(W(\bpi))\subset\bpi\cal Q^-$.
     \item The representations $V(\bpi_{i,a})$, $i\in I$, $a\in\bc^\times$ are minuscule.
      \item For algebras of type $A_n$, $B_n$ and $C_n$ the representations $V(\bpi_{i,a})$ are quasi-minuscule.

      \end{enumerit}
    \end{thm}
The following corollary of the theorem  is immediate from  the
universal property of Weyl modules.
    \begin{cor}\label{cons} Let $V\in\wh{\cal  F}_q(\lie g)$ be an $\ell$--highest weight module with $\ell$--highest weight $\bpi$. Then $$\ch_q(V)=\sum_{\bvpi\in\bpi\cal Q^-}\dim V_\bvpi\bvpi.$$\end{cor}
    The first part of the theorem was proved in \cite{Fre2} and using the Corollary to Theorem 1 they deduced immediately that if $V(\bpi)$ is
     irreducible, then $\ch_q(V(\bpi))=\sum_{\bvpi\in\bpi\cal Q^-}\dim V(\bpi)_\bvpi\bvpi$
      (a stronger general condition on terms of $\ch_q(V(\bpi))$ is proved in
      \cite{her04}). In fact, a stronger result which
implies part (i)   is proved in \cite[Theorem 5.21]{her07}.  In
the case of classical Lie algebras, another proof of part (i) and
a Weyl--character type formula for the $q$--character of the
fundamental representations is given  in \cite{CM}.

Part (ii) of the theorem  is now immediate from   from  Theorem
\ref{presentweyl}.  Part (iii) of the theorem was proved in
\cite{Fre2}. Closed formulae for the characters of fundamental
representations for the classical Lie algebras can be found in
\cite{Fre2, ks, cm2}.  The $q$--character of
fundamental representations of exceptional types were computed in
\cite{Nae, her02, her05}.
    Part (iv) can be established directly from the closed formulae in \cite{ks, cm2}, and
      was first observed (and proved by a different method) in \cite{her05}.
The quasi--minuscule fundamental representations for algebras of
type $D_n$ can be found in \cite{cm2} while the case of $F_4$ and
$G_2$ can be found in  \cite{her02, her05}.

   \subsection{}  We now give  a proof of  parts (i), (iii) and (iv) of the Theorem   for algebras of type $A_n$ since it follows easily from the results discussed so far. The elements $\omega_i\in P^+$, $i\in I$ are minuscule i.e there do no exist dominant integral weights $\mu\in P^+$ such that $\omega_i-\mu\in Q^+\setminus \{0\}$ and $V(\omega_i)$ is a minuscule representation of $\bu_q(\lie g)$. Since $\wt V(\bpi_{i,a})\subset \omega_i-Q^+$, by Theorem \ref{weyldefn}, it follows that $V(\bpi_{i,a})\cong V(\omega_i)$ as $\bu_q(\lie g)$--modules. Hence $$V(\bpi_{i,a})_\mu\ne 0\ \iff \ \mu=w\omega_i,\  \ \
    \dim V(\bpi_{i,a})_{w\omega_i}=1,\ \ w\in W.$$ Using Theorem \ref{braidthm}, we get
    $$\dim V(\bpi_{i,a})_{w\omega_i} =\dim V(\bpi_{i,a})_{\tilde T_w\bpi_{i,a}}, $$ which proves (iv) .
    By Corollary \ref{correctcm}, we see that $\tilde T_w\bpi_{i,a}\in\bpi_{i,a}\cal
    Q^-$.
  Let $W_i=\{w\in W: w\omega_i=\omega_i\}$, and let $W^i$ be
  a set of coset of representatives for $W/W_i$. Then, we have by Theorem
  \ref{braidthm} that
  $$\ch_q V(\bpi_{i,a})=\sum_{w\in W^i} \tilde T_w\bpi_{i,a},$$ and the proof of parts (i) and (iv) of the  theorem is complete.

\subsection{}  Let us consider the case of $V(\bpi_{1,a})$ for $A_n$.
Here, we have $$\wt_\ell
V(\bpi_{1,a})=\{(\bpi_{i,aq^{i+1}})^{-1}\bpi_{i+1, aq^i}: 0\le
i\le n\},$$ where we understand that $\bpi_{n+1,aq^n}=1$. The
following is analogous to Proposition \ref{weylsl2qm} and is
proved in a similar fashion.
 \begin{prop} Suppose that $\lie g$ is of type $A_n$ and that  $\bpi=(\pi,1,\cdots ,1)\in\cal P^+$.
 Then $$\ch_q(W(\bpi))=\sum_{\bpi'\in\cal P^+(\bpi) }d_{\bpi}(\bpi')\bpi',$$ where $$\cal P^+(\bpi)=\left\{(\frac {\pi_1(u)}{\pi_2(q^2u)},\frac{\pi_2(qu)}{\pi_3(q^3u)},\cdots, \frac{\pi_n(q^{n-1}u)}{\pi_{n+1}(q^{n+1}u)}): \pi=\prod_{s=1}^{n+1}\pi_j,\ \pi_j\in\bc[u]\right\},$$
 and  $$d_\bpi(\bpi')= \prod_{s=1}^{k}\prod_{i=1}^{n+1}\binom{r_s-m_{s,1}-m_{s,2}-\cdots -m_{s, i-1}}{m_{s,i}},$$
 where $\pi=\prod_{s=1}^k(1-a_ku)^{r_k}$, $a_k\ne a_m$ if $k\ne m$ and $\pi_j=\prod_{s=1}^k(1-a_ku)^{m_{s,j}}$, $1\le j\le n+1$. In particular, $W(\bpi)$ is quasi--minuscule iff the roots of $\pi$ are distinct.\hfill\qedsymbol
 \end{prop}
\subsection {} Suppose now that $\lie g$ is a simple Lie algebra and assume that $i\in  I$ is such that the coefficient of $\alpha_i$ in the highest root is one.  Then it is known that the corresponding $\bu_q(\lie g)$--representation $V(\omega_i)$ is quasi--minuscule and also that (\cite{Cha2} for instance)  $V(\bpi_{i,a})\cong V(\omega_i)$. In particular, $V(\bpi_{i,a})$ is a quasi--minuscule representation of $\wh{\bu}_q(\lie g)$. 
The converse is also true for type $A_n$ as we see from the preceding discussion : if a fundamental representation is quasi-minuscule then it is simple as a $\bu_q(\lie g)$-module.
For $D_n$ from \cite{cm2} and more generally it has been observed by Nakajima that the converse is also true for all simply--laced algebras.
In the non--simply laced case however, the converse need not be
true, for instance by Theorem \ref{lcone} (iv), all the
fundamental representations of $B_n$ are quasi--minuscule, but
$V(\bpi_{i,a})$ is not irreducible as a $\bu_q(\lie g)$--module if
$i\ne 1,n$, see \cite{Cha2} for instance. Examples for $G_2$ can
be found in \cite{her02}.

\subsection{} As a consequence of Theorem \ref{lcone}(i), more precisely its Corollary \ref{cons}, and by using the grading of $q$-characters by the usual weight lattice of $\lie g$, the $q$-characters of simple representations are clearly $\ZZ$-linearly independent (as for usual characters of simple representations of $\lie g$). Then it follows \cite{Fre} that $\ch_q: \wh{Rep}(\lie g)\to\bz[\cal P]$ is an injective ring homomorphism.

The image of $\ch_q$ has been characterized in \cite{Fre2} as the
intersection of the kernel of screening operators (this is
analogous to the invariance for the Weyl group action of usual
characters). Besides an element of $\text{Im}(\ch_q)$ is characterized \cite{Fre, Fre2} by the multiplicity of the dominant terms $\bpi\in\mathcal{P}^+$.
In other words, we have $\ch_q(V)=\ch_q(W)$ iff $\dim V_\bpi=\dim
W_\bpi$ for all $\bpi\in\cal P^+$.

\section{Affinizations of $\lambda$ and the poset $\bd_\lambda$}

 One tool that has been used very effectively for example in
 \cite{Nae, her02} in the computation of $q$--characters of certain
  families of modules is the Frenkel--Mukhin algorithm \cite{Fre2} and
  we shall describe this in another section. The braid
   group action on $\cal P$ and Theorem \ref{braidthm}
    can be used to compute $q$--characters although this approach
    has not as yet been fully explored outside the fundamental modules.
    As explained in the introduction,
     in the $ADE$-case an algorithm
      is given in \cite{Nab}
      which gives the $q$-character of an
       arbitrary simple finite dimensional representation.
       In practice however, it is difficult to write closed
        formulae for the q-character and
         the first thing  is to identify
          suitable families of modules
          for which these methods can be made to work and possibly lead to   closed formulae.
           The fundamental modules and the Weyl modules introduced
            in the previous section are examples of such families
            and in this section we shall identify some other natural
             families of modules which could provide further examples.
             We shall also be interested in branching rules,
              i.e the decomposition of an object of $\wh{\cal F}_q(\lie g)$ as a direct of sum of simple objects of $\cal F_q(\lie g)$.

   \subsection{}        Since any     $V\in\wh{\cal F}_q(\lie g)$   is completely reducible as a $\bu_q(\lie g)$--module,
   we have $$V\cong_{\bu_q{(\lie g)}} \bigoplus_{\mu\in P^+} V(\mu)^{\oplus {m_\mu}(V)},\ \ m_\mu(V)=
    \dim\Hom_{\bu_q(\lie g)}(V(\mu), V).$$ Set, \begin{equation}\label{chg}\ch_{\lie g}(V)=\sum_{\mu\in P^+}m_\mu(V)\ch(V(\mu)).\end{equation}

  Clearly knowing either $\ch$ or $\ch_q$
 implies that one knows $\ch_{\lie g}$ in principle. In practice
 however it is very hard to find from this information, a closed formula for $\ch_{\lie
 g}$. This is seen in other situations: for instance in the case
 of objects of $\cal F_q(\lie g)$, the character $\ch(V(\lambda))$
 is given by the Weyl character formula and since $\ch$ is
 multiplicative, one knows the character of the tensor product of
 $V(\lambda)\otimes V(\mu)$. But the understanding of $\dim\Hom_{\bu_q(\lie
 g)}(V(\nu), V(\lambda)\otimes V(\mu))$ is a hard problem.
In this  section, we shall see some examples of families of
modules where $\ch_{\lie g}$ is known.
\subsection{}  We shall need two automorphisms introduced in  \cite{Cha0, Chari2} of $\wh{\bu}_q(\lie
g)$. Given $a\in\bc^\times$, there exists an automorphism
$\tau_a:\wh{\bu}_q(\lie g)\to\wh{\bu}_q(\lie g)$ which is defined
on generators by $$\tau_a(x^\pm_{i,r})=a^{\pm r} x^\pm_{i,r},\ \
\tau_a(\psi_{i,\pm m})=a^{\pm r}\psi_{i,\pm m},\ \ \tau_a(k_i^{\pm
1})=k_i^{\pm 1},$$ for all $i\in I$, $r\in\bz$, $m\in\bz_+$.
Let $\sigma: \wh{\bu}_q(\lie
g)\to\wh{\bu}_q(\lie g)$ be the involution satisfying,
$$\sigma(x_{i,r}^{\pm})= x_{i,-r}^{\mp}\text{ , }\sigma(\phi_{i,\pm m}^\pm)= \phi_{i, \mp m}^\mp \text{ , }\sigma(k_i)=k_i^{-1},$$
for $i\in I$, $r\in\bz$, $m\in\bz_+$. Given $V\in\Ob\wh{\cal
F}_q(\lie g)$, let $\tau_a^*(V)$ and $\sigma^*(V)$ be the object
obtained by pulling $V$ back via the automorphism $\tau_a$ and
$\sigma$ respectively.  Define corresponding automorphisms
$\tau_a:\bz[\cal P]\to\bz[\cal P]$ sending  $\bpi\to\bpi_a$ (resp.
$\sigma:\bz[\cal P]\to\bz[\cal P]$ sending $\bpi\to\bpi_\sigma$)
by:
$$\tau_a(\bpi_{i,b})=\bpi_{i,ab},\ \ i\in I,\ \  b\in\bc^\times,$$
$$\sigma(\bpi_{i,b})=\bpi_{-w_\circ(i),q_i^2b^{-1}}, \ \ i\in I,\
\ b\in\bc^\times,$$ where we recall that $w_\circ$ is the longest
element of the Weyl group and that $-w_\circ$ induces an
involution of the Dynkin diagram of $\lie g$. The following can be
found in \cite{Chari2} and \cite[Proposition 5.1]{Cha8}
\begin{prop}\label{isotwist} Let $\bpi\in\cal P^+$ and $a\in\bc^\times$. Then $$\tau_b^*V(\bpi)\cong V(\bpi_b),\ \ \sigma^*(V(\bpi))\cong \tau_h^*V(\bpi_\sigma),$$
for some $h\in \bc^\times$ which is independent of $\bpi$.
\end{prop}
In fact, $h$ is the dual Coxeter number of $\lie g$. The statements of the following corollary were first proved  in \cite{Fre} and  \cite{herma} respectively and are consequence of the fact that for $V\in\wh{\cal F}_q(\lie g)$
 $$\dim (\sigma^* V)_{\sigma\bvpi}= \dim V_\bvpi=\dim (\tau_a^*V)_{\tau_a\bvpi}.\ \ \ \ $$

\begin{cor} Let $\bpi\in\cal P^+$. Then,$${\rm\ch}_q(\tau_a^*V(\bpi))
=\tau_a{\rm\ch}_q(V(\bpi)), \ \ \ \rm{\ch}_q(\sigma^*V(\bpi))=
\sigma({\ch_q}(V(\bpi)).$$
\end{cor}

 \subsection{} Define an equivalence relation $\sim$ on objects of $\wh{\cal F_q}(\lie g)$ by:
  $V\sim V'$ iff $V$ and $V'$ are isomorphic as $\bu_q(\lie g)$--modules and let $[V]$ denote the
   equivalence class of $V$. Note that by Proposition
   \ref{isotwist}, we have $$[V]=[\tau_a^*(V)], \ \
   a\in\bc^\times.$$
It is shown in \cite{Chari2} and can be verified easily, that one
 has a  partial order on the set of equivalence classes  given by: $[V]\ge [V']$ iff for all $\nu\in P^+$, either\\
$\bullet$ $m_\nu(V)\ge m_\nu(V')$, or\\
$\bullet$ there exists $\mu\in P^+$ with $\mu\ge\nu$ such that $m_\mu(V)>m_\mu(V')$.

For $\lambda\in P^+$, set $$\bd_\lambda=\{[V(\bpi)]:
\wt\bpi=\lambda\}.$$ The following is proved in \cite{Chari2}, we
include a proof here, since it follows easily from the results of
the previous sections.
\begin{prop} For all $\lambda\in P^+$ the set $\bd_\lambda$ is  a finite poset.
\end{prop}
\begin{pf}
 Suppose that $\bpi\in\cal P^+$ and $\wt\bpi=\lambda$.
Then $V(\bpi)$ is a quotient of $W(\bpi)$.
Since $W(\bpi)$ is a finite--dimensional
 $\bu_q(\lie g)$--module it follows that $$\ch_{\lie g} W(\bpi)=\sum_{\nu\in P^+} m_\nu(W(\bpi)) V(\nu),$$ and hence we have $$\ch_{\lie g}V(\bpi)=\sum_{\nu\in P^+} m_\nu(V(\bpi)) V(\nu),\ \ m_\nu(V(\bpi))\le n_\nu(V(\bpi)).$$ The proposition follows since there are only finitely many $\nu\in P^+$ with $m_\nu(W(\bpi))\ne 0$.
\end{pf}

\subsection{}  Using Theorem \ref{sl2}, it is easy
to describe the poset $\bd_\lambda$ in the case when $\lie g$ is
of type $\lie{sl_2}$. Suppose that $\lambda=s\omega_1$ and let
$\pi,\pi'$ be such that $\wt\pi=\lambda=\wt\pi'$. Then, it is
trivially seen  that
$$V(\pi)\cong_{\bu_q(\lie{sl}_2)}V(\pi')\ \ \iff
\pi=\prod_{k=1}^K\pi^{m_k}_{a_k},\ \
\pi'=\prod_{k=1}^K\pi^{m_k}_{b_k},$$ where
$\pi=\prod_{k=1}^K\pi^{m_k}_{a_k}$ (resp.
$\pi'=\prod_{k=1}^K\pi^{m_k}_{b_k})$, is the $q$-factorization of
$\pi$ (resp. $\pi'$). In other words, the elements of
$\bd_\lambda$ are just elements of ${\rm{Par}}(s)$, the set
partitions of $s$. Consider the reverse lexicographic ordering on
${\rm{Par}}(s)$, this is the order where $\{s\}$ is the smallest
partition and $\{1,\cdots, 1\}$ the maximal element. It is now an
elementary exercise in the representation theory of $\lie{sl_2}$
to prove the following.
\begin{lem} Assume that $\lie g$ is of type $\lie{sl_2}$ and that $\lambda=s\omega_1$.  Then  $\bd_\lambda$ is isomorphic
as a poset to the  set ${\rm{Par}}(s)$ of partitions of $s$
equipped with the reverse lexicographic order.
\end{lem}

\subsection{} In the case of an arbitrary simple Lie algebra, very little is known about the poset $\bd_\lambda$ except for the minimal and maximal elements of the poset. We first show that the maximal element of $\bd_\lambda$ is just the equivalence class of $W(\bpi)$ for $\bpi\in \cal P^+$ with $\wt\bpi=\lambda$. For this,  it suffices by Theorem \ref{weyldefn} to prove that there exists $\bpi_0\in\cal P^+$ with $\wt(\bpi_0)=\lambda$ such that $W(\bpi_0)$ is irreducible. Write $\lambda=\sum_{i\in I} m_i\omega_i$ and set $\bpi_0=\prod_{i\in I}(\bpi_{i,1})^{m_i}$. By Corollary \ref{fundirr} it follows that $\otimes_{i\in I}V(\bpi_{i,1})^{\otimes m_i}$ is irreducible and by Theorem \ref{presentweyl} we see that it is isomorphic to $W(\bpi_0)$ and we have proved:
\begin{lem} For $\lambda\in P^+$, there exists $\bpi\in\cal P^+$ with $\wt\bpi=\lambda$ such
that $[W(\bpi)]$ is the unique maximal element of $\bd_\lambda$.
\end{lem}

\subsection{} The minimal elements of the poset $\bd_\lambda$ were studied in \cite{Chari2, Cha7, Cha8, Cha9}.
Notice first that if $\lambda=\omega_i$ and $a\in\bc^\times$, then
$\bd_{\omega_i}=\{[V(\bpi_{i,a})]\}$, and in particular any
affinization of $\omega_i$ is minimal. For general $\lambda$, the
picture is more complicated.   Assume  that the
nodes of the Dynkin diagram are numbered as in \cite{Bo}. For
$\lambda\in P^+$, set $$\supp\lambda=\{i\in I: \lambda(i)>0\},$$ and
let $i_{\min\lambda}$ be the minimal element of $\supp\lambda$.

\subsection{}\label{notde}  If  $\lie g$ is not of type $D$ or $E$, then there is a unique minimal element in $\bd_\lambda$. Moreover there exist two elements $\bpi_{\min\lambda}[k]\in\cal P^+$, $k=1,2$ such that:\\
$\bullet$ if $|\supp_\lambda|>1$, there does not exist $a\in\bc^\times$ such that $\bpi_{\min\lambda}[1]=\tau_a(\bpi_{\min\lambda}[2])$,\\
$\bullet$ if $\bpi\in\cal P^+$ is such that $V(\bpi)\sim V(\bpi_{\min\lambda}[k])$ for some $k\in\{1,2\}$, then there exists $a\in\bc^\times$ such that $\bpi=\tau_a(\bpi_{\min\lambda}[k])$.

The elements are defined as  follows:
If $i\notin\supp\lambda$ the $i^{th}$ coordinate of
$\bpi_{\min\lambda}[k]$ is one while if $i\in\supp\lambda$, we let
\begin{equation}\label{min1}
(\bpi_{\min\lambda}[1])_i=\pi_{c_i(\lambda)}^{\lambda(i)},\
 \ c_i(\lambda)=q^{-\sum_{i_{\min} \lambda\leq j\leq i}(d_j\lambda(j)+d_{j+1}\lambda(j+1)+d_{j+1}-a_{j+1,j}-1)},\end{equation}
 \begin{equation}\label{min2}(\bpi_{\min\lambda}[2])_i=\pi_{c_i'(\lambda)}^{\lambda(i)},\ \ \ c_i'(\lambda)=q^{\sum_{i_{\min}\lambda \le j\le i} d_j\lambda(j)+d_{j+1}\lambda(j+1)+d_{j}-a_{j,j+1}-1}.\end{equation}

\subsection{}  We turn now to the case of $D$ or $E$ where the number of minimal
elements in $\bd_\lambda$ can depend on $\lambda\in P^+$. Let $i_0$
be the trivalent node of the Dynkin diagram and let $I_j$, $1\le
j\le 3$ be the disjoint connected components  of
$I\setminus\{i_0\}$.  For $1\le j\le 3$ and $\lambda\in P^+$, set  $\lambda_j=\sum_{i\in I_j}\lambda(j)\omega_j$ and for $\bpi=(\pi_1,\cdots ,\pi_n) \in\cal P^+$ let $\bpi_j$ be defined by setting the $i^{th}$ co-ordinate to be 1 if $i\notin I_j$ and to be $\pi_i$ if $i\in I_j$. Clearly $$\lambda=\lambda(i_0)\omega_{i_0}+\lambda_1+\lambda_2+\lambda_3,\ \ \bpi=\bpi_{i_0}^{\lambda(i_0)}\bpi_1\bpi_2\bpi_3.$$
 Suppose first that $\lambda\in P^+$ is such that
  \begin{equation}\label{antype}  \supp_\lambda\subset I_{j_1}\cup\{i_0\}\cup I_{j_2},\end{equation} for some $j_1,j_2\in\{1,2,3\}.$
   Then the situation is identical to the one discussed in Section \ref{notde}.

 Suppose next that \begin{equation}\label{regular} i_0\in \supp_\lambda\ \ {\rm{and}}\ \  I_j\cap \supp_\lambda\ne\emptyset,\ {\rm{for \ all}}\ \  1\le j\le 3.\end{equation} Then $\bd_\lambda$ has three minimal elements which are obtained as follows.  We have  $[V(\bpi)]$ is a minimal element of $\bd_\lambda$ iff:\\
$\bullet$  $[V(\bpi_j\bpi_{i_0})]$ is a minimal element of $\bd_{\lambda_j+\lambda(i_0)\omega_{i_0}}$ for all $1\le j\le 3$,\\
 $\bullet$ there exists a permutation $\sigma$ of $\{1,2,3\}$ such that  the elements $[V(\bpi_{\sigma(1)}\bpi_{i_0}\bpi_{\sigma(2)})]$ and
  $[V(\bpi_{\sigma(1)}\bpi_{i_0}\bpi_{\sigma(3)})]$ are
   minimal elements of $\bd_{\lambda_{\sigma(1)}+\lambda(i_0)\omega_{i_0}+ \lambda_{\sigma(2)}}$ and
    $\bd_{\lambda_{\sigma(1)}+\lambda(i_0)\omega_{i_0}+ \lambda_{\sigma(3))}}$ respectively.\\
      Since $\lambda_{j}+\lambda(i_0)\omega_{i_0}$ and $\lambda_{j}+\lambda(i_0)\omega_{i_0}+\lambda_{k}$
      for $1\le j\ne k\le 3$ satisfy the condition
      in \eqref{antype} one can write down the elements $\bpi$ explicitly  using \eqref{min1} and \eqref{min2} and we find that there are three minimal elements.

The remaining case when $i_0\notin \supp_\lambda$ has been studied
when $\lie g$ is of type $D_4$ in \cite{Cha9}. Here, one finds
that the number of minimal elements in $\bd_\lambda$ \lq\lq
increases\rq\rq with $\lambda$.  Virtually nothing is known in
this case beyond $D_4$.

\subsection{}\label{dkr} Consider the special case $\lambda=m\omega_i$.
 The preceding results show that the poset $\bd_{m\omega_i}$ has a
    unique minimal element and that $\bpi_{\min{m\omega_i}}=(\pi_1,\cdots,\pi_n)$
     where $$\pi_j(u)=\begin{cases} 1,\ \ \text{ if }j\ne i,\\ 
     \prod_{k=1}^m(1-q_i^{m-2k+1}u)\text{ if }j = i.\end{cases}$$
\begin{defn} The modules $V(\tau_a^*(\bpi_{\min{m\omega_i}}))$, $a\in\bc^\times$ are called the
      Kirillov--Reshetikhin modules.
 \end{defn}
We shall discuss these further  in the last section.

      \subsection{}\label{avant} We now  discuss the
      $\ch_{\lie g}$--characters of minimal affinizations.  In the
      case when $\lie g$ is of type $\lie{sl}_{n+1}$, one can prove
      \cite{cpsmall}
      by using the evaluation homomorphism defined by
    Drinfeld and Jimbo that
      \begin{equation}\label{aniso} V(\bpi_{\min\lambda}[k])\cong_{\lie{sl}_{n+1}} V(\lambda),\ \
      \lambda\in P^+.\end{equation} In the case when $\lie g$ is of type $C_2$,
      it was proved in \cite{Chari2} that
      $$V(\bpi_{\min\lambda}[k])\cong_{\lie g}
      \bigoplus_{r=0}^{s}
      V(\lambda-2r\lambda(2)\omega_2),$$ where $s$ is the integer part of
      $\lambda(2)/2$.
For other Lie algebras, virtually nothing is known about the
   $\ch_{\lie g}$ --character of the minimal affinizations in general.

    In the case
     when
    $\lambda=m\omega_i$, and  $i$ is of classical type, i.e, $\alpha_i$ occurs in the highest root $\theta$ with
     multiplicity at most two, the $\ch_{\lie g}$ --character is studied in \cite{c0}  and the results confirm the conjectures in  \cite{hkoty, hkott, kl, kr}.
      As an example we give the $\lie g$--structure of in the case when $i=2$ for the algebras of type $B_n$: here we have that
     $$V(\bpi_{\min m\omega_2}) \cong_{\lie g} \bigoplus_{r=0}^{m} V(r\omega_2).$$ Graded versions and classical analogs of these results have been studied in \cite{cm3, cm4}. Other approaches can also be found in \cite{cg}. If $i\in I$ is not classical, then the corresponding Lie algebra is exceptional
     and such a decomposition for $\ch_{\lie g}(V(\bpi_{\min\omega_i}))$ is not known. In the case
     of $E_8$ for instance, even conjectural decomposition formulas are not known
     when $i$ is the trivalent node and $m>1$. However,  character formulas in a different form have been proved (see the last section).

\section{Properties of  minimal affinizations}

In this section,  we  review results which give sufficient
conditions for minimal affinizations to be prime, minuscule or
quasi--minuscule.

\subsection{}

We can prove the following.
\begin{prop} \label{prime} Let $\lambda\in P^+$  and suppose that $a_{ij}\in\{0,-1\}$ for all
$i,j\in\supp\lambda$ and $i\ne j$. Then $V(\bpi_{\min\lambda})$ is
prime.
\end{prop}

\begin{pf} Notice first that under the hypothesis on $\lambda$, the elements $\bpi_{\min\lambda}[k]$, $k=1,2$ are given
explicitly in  \eqref{min1},\eqref{min2}. If $\lie g$ is of type
$A_n$, then one knows from \eqref{aniso} that
$V(\bpi_{\min\lambda}[k])\cong V(\lambda)$ as $\bu_q(\lie
g)$--modules. It is an elementary exercise to see that
$V(\lambda)$ can never be written as a tensor product of two
non--trivial representations of $\bu_q(\lie g)$ and hence the
proposition follows in this case. For the general case, suppose
that there exists $\bpi,\bpi'\in\cal P^+$, such that
$V(\bpi_{\min\lambda}[k])\cong V(\bpi)\otimes V(\bpi')$. Let $I_0$ be
the connected component of $I$ containing $\supp\lambda$. The
subalgebra $\wh{\bu}_0$ generated by $\wh{\bu}_i$, $i\in I_0$ is
isomorphic to $\wh{\bu}_q(\lie sl_{r+1})$, where $r=|I_0|$. Using
the formulae for comultiplication given in Lemma \ref{comult}, one
sees that
$$\wh{\bu_0}(v_{\bpi_{\min\lambda}[k]})\cong\wh{\bu}_0v_\bpi\otimes\wh{\bu}_0
v_{\bpi'},$$ which means that $\wh{\bu}_0v_{\bpi_{\min\lambda[k]}}$ is not a minimal
affinization for $\lambda$, where we regard $\lambda$ as  an
element of $P_0^+$ (the weight lattice corresponding to $\lie{
sl}_{r+1}$). But again, an inspection of \eqref{min1},
\eqref{min2} shows that this is a contradiction.
\end{pf}
\begin{rem} More generally, one can prove that any minimal
affinization is prime and  a proof will appear elsewhere. The
converse statement however is not true, examples of prime
representations which are not minimal can be found for $\lie sl_3$
in \cite{cpp}.\end{rem}

\subsection{} We explain two essential tools that are needed to continue our study. The first is
an
algorithm defined   in  \cite{Fre2} which can be used to compute
the $q$-character of minuscule representations.

Given
$\bpi=(\pi_1,\cdots,\pi_n)\in\cal P^+$,
 let  $\cal Q_\bpi^-$ be the
submonoid of $\cal Q^-$ generated by the elements
$$\{(\balpha_{i,aq^m})^{-1}: i\in I,\  m\in\bz, \ a\in\bc^\times,\ \prod_{i=1}^n\pi_i(a^{-1})=0\}.$$
The set $\bpi\cal Q^-_\bpi$ is countable and we  fix an
enumeration of this set $\{\bvpi_r\}_{r\ge 0}$ so that,\\
$\bullet$ $\bvpi_0=\bpi$,\\
$\bullet$  $r\ge r'$ implies
$\bvpi_r\preccurlyeq\bvpi_{r'}$.

For $i\in I$, define $\bop_i: \cal P\to\bz[\cal P]$ by: for
$\bvpi=(\varpi_1,\cdots,\varpi_n)$, we have
$$\bop_i(\bvpi)=\begin{cases} 0\ \ \varpi_i\notin\bc[u],\\
(\varpi_1,\cdots,\varpi_{i-1},
1,\varpi_{i+1},\cdots,\varpi_n)\ch_q^i V(\varpi_i)\ \
\varpi_i\in\bc[u],\end{cases}$$ where $\ch_q^i$ is the $q$--character of the module $V(\varpi_i)$ for $\wh{\bu}_i$,
expect that in $(\varpi_i)^{-1}\ch_q^i$ we use the $\balpha_{i,a}^{-1}$ instead of the roots $\balpha_a^{-1}$ of $\wh{\bu}_i$.
If $\bvpi'\in\cal P$ we
let $\bop_i(\bvpi)_{\bvpi'}$ be the coefficient of $\bvpi'$ in
$\bop_i(\bvpi)$. Note that $\bop_i(\bvpi)_{\bvpi'}\in\bz_+$.

For $i\in I$ and $r\in\bz_+$ define integers $s(\bvpi_r) $,
$s_i(\bvpi_r)$, inductively and simultaneously by:

 \begin{gather*} s(\bvpi_0)=1,\ \
s_i(\bvpi_0)=0,\\ s_i(\bvpi_r)=\underset{r' <
r}{\sum}(s(\bvpi_{r'})-s_i(\bvpi_{r'}))[\bop_i(\bvpi_{r'})]_{\bvpi_r}\text{
 },\ \ r\ge 1, \\  s(\bvpi_r) = \text{Max}_{j\in
I}(s_j(\bvpi_r)),\ \ \ r\ge 1\end{gather*} Finally, set $FM(\bpi)
= \sum_{r\geq 0}s(\bvpi_r) \bvpi_r$.

\begin{defn} The Frenkel-Mukhin algorithm is said to be well-defined if $FM(\bpi)\in\bz[\cal
P]$, i.e $s(\bvpi_r)=0$ for all but finitely many $r$. We say that
$V(\bpi)$ is regular if $\ch_q(V(\bpi)) = FM(\bpi)$.
\end{defn}

\begin{thm}\cite{Fre2} A minuscule simple module is regular.
\end{thm}
In general a simple representation is not regular (examples were
first considered in \cite{nn4} for type $C_3$). Note that if
$V(\bpi)$ is regular, then a representation theoretical
interpretation of the integers $s(\bpi_r) - s_i(\bpi_r)$ can be
found in
 \cite{small}.

In \cite{her02}, a $F(\bpi)\in \bz[\cal
P]$ has been constructed 
for any dominant $\bpi$.
If $V(\bvpi)$ is minuscule, then $F(\bvpi)=FM(\bvpi)=\ch_q(V(\bvpi))$. For general $\bvpi$,
$F(\bvpi)\in\text{Im}(\ch_q)$, but may have
negative coefficients. $FM(\bvpi)$ has nonnegative coefficients, but may not
belong to the image of $\ch_q$.

\subsection{}

The Frenkel-Mukhin algorithm produces elements of $\cal P$ which
could occur in the set $\wt_\ell V(\bpi)$.  We now give a result
which gives a sufficient condition \cite{small} for an element of
$\cal P$ to not be in the set $\wt_\ell V(\bpi)$. This elimination theorem is
useful to prove that a representation is minuscule. Another
application to the smallness conjecture is given later.

Define  a morphism of monoids $\Ht:\cal Q^\pm\to\bz_+$ by
extending
$$\Ht(\balpha_{i,a})^{\pm 1}=1,\ \ i\in I, a\in\bc^\times.$$ Let $\cal
P_i^+=\{\bvpi\in\cal P: \pi_i\in\bc[u]\}.$ The following is proved
in \cite{small}.
\begin{thm}\label{racourc} Let $\bpi\in\cal P^+$ and assume that
$\bvpi\preccurlyeq \bpi$ satisfies the following conditions for
some $i\in I$.
\begin{enumerit}
\item [\rm(i)] There exists a unique element $\tilde\bvpi\in
\wt_\ell(V(\bpi))\cap \cal P_i^+\cap\bvpi\cal Q^-$ such that
$\tilde\bvpi\ne\bvpi$, and its multiplicity is $1$
, \item[\rm(ii)] $ x_{i,r}^+
(V(\bpi)_{\tilde\bvpi})=\{0\}$ for all $r\in\bz$, \item[\rm (iii)]
$\bop_i(\tilde\bvpi)_{\bvpi}=0$, \item[\rm (iv)] if $\bvpi'$ is
such that $V(\bpi)_{\bvpi'}\cap\wh{\bu}_i V_{\tilde\bvpi}\ne 0$
 and $\bvpi_i'\in\bc[u]$, then
 $\Ht(\bvpi'\bpi^{-1}) \geq \Ht(\bvpi\bpi^{-1})$,
\item[\rm (v)] for all $j\neq i$, we have
$$\{\bvpi'\in \wt_\ell(V(\bpi)): \Ht(\bvpi'\bpi^{-1}) < \Ht(\bvpi\bpi^{-1})\}\cap
\bvpi\cal Q_j=\emptyset,$$ where $\cal Q_j$ is the subgroup of
$\cal Q$ generated by the elements $\balpha_{j,a}$,
$a\in\bc^\times$. \end{enumerit}Then $V(\bpi)_\bvpi=0$.\end{thm}

\subsection{}\label{met}
We now consider the problem of giving  explicit formulae for the $q$--character of a minimal
affinization. It is useful to see the introduction of \cite{herma}
for references on previous results in this direction; in
particular in type $A$ the results can be found in \cite{Che1, nt}, in the case of Yangians and in  \cite{Fre3} for quantum
affine algebras. In the special case of the  Kirillov-Reshetikhin,  the results can be extracted from \cite{Nad, her06}
and we refer to the next  section, in which they are discussed in greater detail, for references.

 The
following result is proved in \cite{herma}.
\begin{thm} \label{elim}Let $\lambda\in P^+$ and assume that $\bpi_{\min\lambda} \in\{\bpi_{\min\lambda}[k]: k=1,2\}$.
\begin{enumerit}
\item If $\lie g$ is of type $A_n$, $B_n$ or $G_2$, then
$V(\bpi_{\min\lambda})$ and $\sigma^*V(\bpi_{\min\lambda})$ are
minuscule. \item if $\lie g$ is of type $C_n$ and $\lambda(n)=0$,
then $V(\bpi_{\min\lambda})$ is minuscule if $\bpi_{\min\lambda}$
satisfies \eqref{min1} and  $\sigma^*V(\bpi_{\min\lambda})$ is
minuscule if $\bpi_{\min\lambda}$ satisfies \eqref{min2}. An
analogous result holds if $\lie g$ is of type $F_4$ if we assume
that $\lambda(4)=0$. \item If $\lie g$ is of type $D_n$ and 
$\lambda(n-1)=\lambda(n)$, then
$V(\bpi_{\min\lambda})$ is minuscule if $\bpi_{\min\lambda}$
satifies \eqref{min1} and $\sigma^*V(\bpi_{\min\lambda})$ is
minuscule if $\bpi_{\min\lambda}$ satisfies \eqref{min2}.
\hfill\qedsymbol
\end{enumerit}
\end{thm}

The main points of the  proof are a generalization of the methods of
\cite{her06}. There are two main steps: the first one is to
compute the terms which occur ``at the top'' of the
$q$-characters, that is to say the first few terms. This is done by Theorem \ref{racourc}
 to show that some terms cannot occur in the
$q$-character. Of the terms that occur, only the highest weight is in $\cal P^+$.
 Then most have the following right-negative property :

\begin{defn}\label{monomrn}\cite{Fre2} A non trivial $\bvpi = \prod_{i\in I,a\in\CC^*}\bomega^{u_{i,a}(\bvpi)}$ is said to be right-negative if for all $a\in\CC^*, j\in I$ we have $u_{j,aq^{L_a}}(\bvpi)\neq 0\Rightarrow u_{j,aq^{L_a}}(m)<0$ where 
$$L_a=\text{max}\{l\in\ZZ/\exists i\in I, u_{i,aq^L}(\bvpi)\neq 0\}.$$\end{defn}

The second step is to use the information
of the top of the $q$-character, the right-negative property, to prove that all other elements in $\wt_\ell(V(\bpi))$ also have the right negative property. This implies immediately, that they are not elements of $\cal P^+$. For this step, we need  a
representation theoretical induction inside the module and a
 crucial ingredient is the structure of Weyl-module in the $sl_2$-case.

Although, the theorem is not the best possible in the case of
algebras of type $C_n$, $D_n$ and $F_4$, it is  not true that all
minimal affinizations are minuscule.  For instance, when $\Glie$
of type $C_3$, and $\lambda=2\omega_2+\omega_3$, then one can see
that a  corresponding minimal affinization is not minuscule
although it satisfies \eqref{min1}. Other counter-examples may be
found in \cite{herma}.

\subsection{} Recall the classical result
that the Jacobi--Trudi determinant gives the character of the
irreducible representation $V(\lambda)$ of $\lie{sl}_{n+1}$. A
generalization of the Jacobi--Trudi determinant to other classical
Lie algebras in terms of tableaux can be found in \cite{kos} for
type $B$, and \cite{nn1, nn2, nn3} for general
classical type. In \cite[Conjecture 2.2]{nn1} Nakai and Nakanishi conjectured that when $\lie g$ is of  classical type,
the Jacobi-Trudi type determinant gives  the $q$-character of an
irreducible representation of the corresponding   quantum affine
algebra.  The following is proved in \cite{herma} and confirms
their conjecture when $\lie g$ is of type $A$ or $B$.

\begin{thm}\label{nnconj} Assume that $\Glie$ of type $A$ or $B$ and let $\lambda\in P^+$ be such that $\lambda(n)$ is even if $\lie g$ is of type $B_n$.
 The $q$-character of the  minimal affinization  of $\lambda$  is given by the corresponding
Jacobi-Trudi determinant.  In particular, the corresponding
minimal affinizations are
quasi-minuscule.\hfill\qedsymbol\end{thm}

 If $\Glie$ of type $C_4$, for instance
 the minimal affinization of $2\omega_3$ is not quasi--minuscule.
 As discussed earlier in the paper, for $D_n$, there exist
 fundamental representations which are not quasi--minuscule.

\subsection{}   As an illustration, we  give the
q-characters   for type $B_n$ predicted by the conjecture of
\cite{nn2}.

 Recall that a partition $\lambda =
(\lambda_1,\lambda_2,\cdots)$ is a sequence of weakly decreasing
non-negative integers with finitely many non-zero terms. The
conjugate partition is denoted by $\lambda' =
(\lambda_1',\lambda_2',\cdots)$. Given two partitions $\lambda$
and $\mu$, we say that $\mu\subset \lambda$ if $\lambda_i\geq
\mu_i$ for all $i\ge 1$ and  the corresponding skew diagram
denoted $\lambda/\mu$ is
$$\lambda/\mu = \{(i,j)\in\bn\times\bn: \mu_i + 1\leq j\leq \lambda_i\} = \{(i,j)\in\bn\times\bn: \mu_j' + 1\leq i\leq
\lambda_j'\}, $$ and let $d(\lambda/\mu)$ be the length of the
longest column of $\lambda/\mu$.

Assume now that  that $d(\lambda/\mu) \leq n$  and also that
$\lambda/\mu$ is connected (i.e. $\mu_i + 1\leq \lambda_{i+1}$ if
$\lambda_{i+1}\neq 0$).

Let $\mathbf B = \{ 1,\cdots,n,0,\overline{n},\cdots,
\overline{1}\}$. We give the ordering $\prec$ on the set $\mathbf B$
by
\begin{equation*}
  1 \prec 2 \prec \cdots \prec n \prec 0 \prec \overline{n} \prec\cdots \prec \overline{2}\prec \overline{1}.
\end{equation*}
As it is a total ordering, we can define the corresponding maps
$\text{succ}$ and $\text{prec}$. For $a\in\CC^*$, let (we use the notations of \cite{Fre} explained in Section \ref{rappel})
$$ \ffbox{i}_a = Y_{i-1,aq^{2i}}^{-1} Y_{i,aq^{2(i-1)}}\text{ , }\ffbox{\overline{i}}_a = Y_{i-1,aq^{4n - 2i - 2}} Y_{i,aq^{4n - 2i}}^{-1}\text{ for $2 \le i \le n-1$,}$$
$$\ffbox{n}_a
     = Y_{n-1,aq^{2n}}^{-1} Y_{n,aq^{2n - 1}}Y_{n,aq^{2n - 3}}\text{ , }
     \ffbox{\overline{n}}_a = Y_{n-1,aq^{2n - 2}} Y_{n,aq^{2n+1}}^{-1}Y_{n,aq^{2n-1}}^{-1},$$
$$\ffbox{0}_a = Y_{n,aq^{2n+1}}^{-1} Y_{n,aq^{2n-3}},$$
(we denote $Y_{0,a} = Y_{n+1,a} = 1$). For $T =
(T_{i,j})_{(i,j)\in\lambda/\mu}$ a tableaux of shape $\lambda/\mu$
with coefficients in $\mathbf B$, let
$$m_{T,a} = \prod_{(i,j)\in\lambda/\mu} \ffbox{T_{i,j}}_{aq^{4(j-i)}}\in\Yim.$$
Let $\text{Tab}(B_n,\lambda/\mu)$ be the set of tableaux of shape
$\lambda/\mu$ with coefficients in $\mathbf B$ satisfying :
$$[T_{i,j}\preceq T_{i,j+1}\text{ and }(T_{i,j},T_{i,j+1})\neq (0,0)]\text{ and }[T_{i,j}\prec T_{i+1,j}\text{ or }(T_{i,j},T_{i+1,j}) = (0,0)].$$
The tableaux expression of the Jacobi-Trudi determinant \cite{kos,
nn1} is :
$$\chi_{\lambda/\mu,a} = \sum_{T\in \text{Tab}(B_n,\lambda/\mu)} m_{T,a} \in \Yim .$$

Note that we get minimal affinizations of Theorem \ref{nnconj} for $\mu = 0$. The highest weight is $\sum_j \sum_{\{i|\lambda_i = j\}} (3 - r_j)\omega_j$. The proof of the conjecture for more general representations
will appear elsewhere.

\section{Kirillov-Reshetikhin modules}

We conclude this paper
 with a discussion of the Kirillov--Reshetikhin (KR) modules. These were first
 introduced in \cite{kr} and since then have been widely studied.
  They have a number of interesting properties, some of which we have already seen.
  They are the minimal affinizations $W_{m,a}^{(i)} = \tau_a^*V(\bpi_{\min m\omega_i})$ of $m\omega_i$, $m\in\bz_+$, $i\in I$, $a\in\bc^\times$ of section \ref{dkr} and are prime by Proposition \ref{prime}.

  We have also seen  in Theorem \ref{elim} that in
     some cases they are minuscule and we have seen in section \ref{avant} that closed formulas
     are known for $\ch_{\lie g}$ when $i$ is of classical type.
     Some of these results can be improved for these modules as
     we shall see below. In
particular it is one of the first infinite family of simple finite
dimensional representations of $\wh{\bu}_q(\lie g)$ where explicit
uniform character formulas can be given for all types : this important
property of KR modules is the KR conjecture proved by Nakajima (ADE case)
and the second author (general case) and discussed in this section.

     Another important property that
      is expected is that the modules $V(\bpi_{\min m\omega_i})$
        have a crystal basis (for a choice of the spectral parameter) and moreover that a module $V(\bpi)$ has a crystal basis if and only if it is a tensor product of
modules of the form $V(\bpi_{\min m\omega_i})$.  Another
 important motivation for the study of the  Kirillov--Reshetikhin modules
  is their connections with solvable lattice models \cite{hkoty, hkott}. There is extensive literature on the subject and we just give a quick overview
  in this section, with pointers to the appropriate references.

\subsection{} By Proposition \ref{isotwist} we see that $W_{m,a}^{(i)}\cong_{\lie g}
W_{m,1}^{(i)}$ and hence $\ch W_{m,a}^{(i)}=\ch W_{m,1}^{(i)}$. The first part of the Kirillov-Reshetikhin conjecture gives a closed formula for the character $\ch$ of an arbitrary tensor product of the modules $W_{m,1}^{(i)}$.

   For a sequence $\nu=(\nu_k^{(i)})_{i\in
I, k> 0}$ of non--negative integers,  such that  all but finitely many $\nu_k^{(i)}$ are  zero we set :
$$\mathcal{F}(\nu)=\underset{N=(N_k^{(i)})}{\sum}\underset{i\in I, k>0
}{\prod}\begin{pmatrix}P_k^{(i)}(\nu,N)+N_k^{(i)}\\N_k^{(i)}\end{pmatrix}e(-k
N_k^{(i)}\alpha_i)$$ where
$$P_k^{(i)}(\nu,N)=\underset{l=1...\infty}{\sum}\nu_l^{(i)}\text{min}(k,l)-\underset{j\in
I,l>0}{\sum}N_l^{(j)}r_iC_{i,j}\text{min}(k/r_j,l/r_i),$$
$$\begin{pmatrix}a\\b\end{pmatrix}=\frac{\Gamma(a+1)}{\Gamma(a-b+1)\Gamma(b+1)},$$ and $\Gamma$ is the usual gamma function.
\begin{thm}[The KR conjecture]\label{conjun} For a sequence $\nu=(\nu_k^{(i)})_{i\in I, k> 0}$ such that for all but finitely many
$\nu_k^{(i)}$ are zero,  we have : \begin{equation}\label{cha} \underset{i\in I, k\geq
1}{\prod}(\ch(W_{k,1}^{(i)}))^{\nu_k^{(i)}}\underset{\alpha\in\Delta_+}{\prod}(1-e(-\alpha))=\mathcal{F}(\nu).\end{equation}\hfill\qedsymbol\end{thm}
The following theorems \ref{tsyst}, \ref{formerkr}, \ref{conv} (that imply Theorem \ref{conjun} by \cite{hkoty, hkott, knt}) are due to  Nakajima \cite{Nab, Nad} for simply--laced algebras (the proof uses geometric methods) and in  full generality to the  second author \cite{her06} (the general proof uses  purely algebraic different methods described in section \ref{met}).

In both proofs the crucial step is the following :
\begin{thm}\label{formerkr} The  modules $W_{k,a}^{(i)}$ are minuscule.\end{thm}

Recently it was proved combinatorially in \cite{dk} that \eqref{cha} can be rewritten in a different form with positive coefficients.

\subsection{} We now discuss the relationship of Theorem \ref{conjun} to $T$ and $Q$ systems.
The $T$--systems were originally introduced in \cite{kns} as functional relations. They can also be viewed as a system of induction relations on
the characters of the Kirillov--Reshetikhin modules and to do this,  we introduce the following representations. For $i\in I$, $k\geq 1$,
$a\in\CC^*$ define the $\wh{\bu}_q(\lie g)$-module $S_{k,a}^{(i)}$ by :
\begin{equation*}
\begin{split}S_{k,a}^{(i)}= \begin{cases}(\underset{j/a_{j,i}=-1}{\bigotimes}W_{k,aq_i}^{(j)})\otimes(\underset{j/a_{j,i}\leq
-2}{\bigotimes}W_{d_ik,aq}^{(j)})&\text{ if $d_i\geq 2$,}
\\(\underset{j/a_{i,j}=-1}{\bigotimes}W_{k,aq}^{(j)})\otimes
(\underset{j/a_{i,j}=-2}{\bigotimes}W_{r,aq}^{(j)}\otimes W_{r,aq^3}^{(j)})&\text{ if $d_i = 1$, $\lie g\neq G_2$, $k = 2r$,}
\\(\underset{j/a_{i,j}=-1}{\bigotimes}W_{k,aq}^{(j)})\otimes
(\underset{j/a_{i,j}=-2}{\bigotimes}W_{r+1,aq}^{(j)}\otimes W_{r,aq^3}^{(j)})&\text{ if $d_i = 1$, $\lie g\neq G_2$, $k = 2r
+1$,}\end{cases}\end{split}
\end{equation*}
and for $r_i = 1$ and $\lie g$ of type $G_2$ ($j\neq i$ is the other node) :
\begin{equation*}
\begin{split}
S_{k,a}^{(i)}=\begin{cases}W_{r,aq}^{(j)}\otimes W_{r,aq^3}^{(j)}\otimes W_{r,aq^5}^{(j)}&\text{ if $k = 3r$,}
\\W_{r+1,aq}^{(j)}\otimes W_{r,aq^3}^{(j)}\otimes W_{r,aq^5}^{(j)}&\text{ if $k = 3r + 1$,}
\\W_{r+1,aq}^{(j)}\otimes W_{r+1,aq^3}^{(j)}\otimes W_{r,aq^5}^{(j)}&\text{ if $k = 3r + 2$.}
\end{cases}
\end{split}
\end{equation*}

\noindent The $S_{k,a}^{(i)}$ are well-defined as the modules in
the definition commute for $\otimes$. Moreover $S_{k,a}^{(i)}$ is
minuscule and so is simple. We denote by $[V]$ the image in
$\wh{\text{Rep}}(\lie g)$ of a module $V$.

\begin{thm}[The $T$-system]\label{tsyst} Let $a\in\CC^*$, $k\geq 1$, $i\in I$. We have in $\wh{\text{Rep}}(\lie g)$:
$$[W_{k,a}^{(i)}][W_{k,aq_i^2}^{(i)}] = [W_{k+1,a}^{(i)}][W_{k-1,aq_i^2}^{(i)}]
+ [S_{k,a}^{(i)}].$$ \end{thm} The $T$-system holds in the Grothendieck ring, but it can also be written as an exact sequence of
representations.
 The $T$-system implies the $Q$-system which is just the corresponding statement in $\text{Rep}(\lie g)$.

For example in the case of $sl_2$, the $T$-system is just the
following relation which can be easily checked by using the
explicit formulas given above :
$$[W_{k,a}][W_{k,aq^2}] = [W_{k+1,a}][W_{k-1,aq^2}] + 1.$$
The $Q$-system is
$$Q_k^2 = Q_{k+1}Q_{k-1} + 1$$
where $Q_k = \text{ch}(W_{k,a})$ does not depend on the spectral parameter $a$. This just an elementary relation between the characters $\ch(Q_k) = e(k\omega) +
e((k-2)\omega) + \cdots + e(-k\omega)$.

\subsection{} A convergence property of the $q$-characters of KR modules holds :

\begin{thm}\label{conv} The normalized $q$-character of $W_{k,a}^{(i)}$ considered as a polynomial in $\balpha_{j,b}^{-1}$ has a limit as a
formal power series :
$$\exists \underset{k\rightarrow
\infty}{\text{lim}}\frac{\text{ch}_q(W_{k,aq_i^{-2k}}^{(i)})}{\bpi_{k,aq_i^{-2k}}^{(i)}}\in\ZZ[[\balpha_{j,aq^m}^{-1}]]_{j\in I, m\in\ZZ},$$
where $\bpi_{k,aq_i^{-2k}}^{(i)}$ is the highest term of $W_{k,aq_i^{-2k}}^{(i)}$.
\end{thm}
As a direct consequence, a convergence property holds for the characters of KR module : $\mathcal{Q}_k^{(i)}=e(-k\omega_i)\ch(W_{k,a}^{(i)})$ considered as a polynomial in $e(-\alpha_j)$ has a limit as a formal
power series :
$$\exists \underset{k\rightarrow \infty}{\text{lim}}\mathcal{Q}_k^{(i)}\in\ZZ[[e(-\alpha_j)]]_{j\in I}.$$
With the $Q$-system, these representation theoretical results imply by
combinatorial arguments \cite{hkoty, hkott, knt} explicit character in Theorem \ref{conjun}.
\cite{dk} deals with the problem of rewriting (\ref{cha}) into a different expression.

\subsection{} It is expected that KR modules (for a special choice of the spectral parameter) have a crystal basis. This is known for fundamental
representations \cite{kas} (see \cite{hn} and references therein for explicit descriptions). As an application of the branching rules of KR
modules discussed in Section \ref{avant} (the branching rules in \cite{her09} for twisted cases), the conjecture about crystal basis has been proved for classical types (see \cite{os} and references therein).

\subsection{} Let us now go to the question of $q$-characters of KR modules. The fact that they are minuscule implies that they are regular and that their $q$--character can in principle be calculated by using the  Frenkel-Mukhin algorithm. In classical
types, there are  explicit formulas in \cite{kos, knh} (the formulas for fundamental representations are given in \cite{ks, cm2}) which follow from the minuscule property. But explicit formulas for their $q$-character are not known in other cases. It would be interesting as well, to give analogs of fermionic formulas
for their $q$-characters.

\subsection{}
In simply-laced cases, Nakajima \cite{Nab} defined $t$-analogs of
$q$-characters (see \cite{her02} for non simply laced cases based
on a different proof of the existence). Nakajima's construction of
$q,t$-characters is closely related to the geometry of quiver
varieties. The geometric small property (Borho-MacPherson) of
projective morphisms implies a description of their singularities
in terms of intersection homology. This notion for certain
resolutions of quiver varieties \cite{Nab} (analogs of the
Springer resolution) can be translated in terms of
$q,t$-characters. Then by using a modification of the proof of
Theorem in \cite{Nab}, it is proved in \cite{small} that we have
the following purely representation theoretical characterization
of small modules involving $q$-characters without
$q,t$-characters. We will use it as a definition :

\begin{thm}\label{repthchar} Let $\bpi\in\mathcal{P}^+$. $V(\bpi)$ is small if and only if for all $\bpi'\in\mathcal{P}^+$ satisfying $\bpi'\preceq \bpi$, $V(\bpi')$ is minuscule.\end{thm}

Note that a small module is necessarily minuscule. From the
geometric point of view it is important to determine which modules
are small. In particular, Nakajima \cite[Conjecture 10.4]{Nab}
raised the problem of characterizing the Drinfeld polynomials of
small standard modules corresponding to KR modules. The main
result of \cite{small} is an explicit answer to this question
(Theorem \ref{mainres}). First let us note in general the standard
modules corresponding to KR modules are not necessarily small :
for type $A_3$, $V(\bpi_{2,1}\bpi_{2,q^2}\bpi_{2,q^4})$
is not small since  $V(\bpi_{1,q}\bpi_{3,q}\bpi_{2,q^4})$ is not minuscule and $\bpi_{1,q}\bpi_{3,q}\bpi_{2,q^4}\preccurlyeq\bpi_{2,1}\bpi_{2,q^2}\bpi_{2,q^4}$.

\subsection{} Let us give a characterization of small KR modules.
\begin{defn} A node $i\in \{1,\cdots ,n\}$ is said to be extremal (resp. trivalent) if there is a unique $j\in I$ (resp. three distinct $j,k,l\in I$) such that $a_{i,j} < 0$ (resp. $a_{i,j}<0$, $a_{i,k}<0$ and $a_{i,l}<0$).

For $i\in I$, we denote by $k_i$ the minimal $k\geq 1$ such that
there are distinct $i = i_1,\cdots,i_k\in I$ satisfying
$a_{i_j,i_{j+1}} < 0$ and $i_k$ is trivalent.
If there does not exist such $k$, set $k_i = +\infty$.
\end{defn}
For example for $\Glie$ of type $A$, we have $k_i = +\infty$ for
all $i\in I$.

 \begin{thm}\label{mainres}[Smallness problem]\cite{small} Let $k\geq 0, i\in I, a\in\CC^*$. Then $W_{k,a}^{(i)}$ is small if and only if $k\leq 2$ or ($i$ is extremal and $k\leq k_i+1$).\end{thm}

In particular for $\Glie = sl_2$ or $\Glie = sl_3$, all KR modules
are small (it proves the corresponding \cite[Conjecture
10.4]{Nab}). In general it gives an explicit criterion so that the
geometric smallness holds.
A few words about the proof : the \lq\lq only if" part is proved by writing down explicitly an element $\bpi\in\cal P^+$ so that $\bpi\preccurlyeq \tau_a\bpi_{\min k\omega_i}$ which proves that the module is not minuscule. For the \lq\lq if" part, all dominant monomials lower than $\tau_a\bpi_{\min k\omega_i}$ are computed, and then it is proved by using the elimination theorem that they correspond to minuscule representations.

\end{document}